\documentclass[11pt,twoside]{article}
\usepackage{latexsym}
\usepackage{amsfonts}
\usepackage{amsmath}
\usepackage{amsthm}
\usepackage{amssymb}
\usepackage[all]{xy}
\usepackage{graphicx}

\newtheorem{theorem}{Theorem}
\newtheorem{proposition}[theorem]{Proposition}
\newtheorem{definition}[theorem]{Definition}
\newtheorem{example}[theorem]{Example}
\newtheorem{remark}[theorem]{Remark}

\newcommand{\G}{{\mathfrak g}}

\newcommand{\D}{{\cal D}}
\newcommand{\A}{{\cal A}}
\newcommand{\V}{{\cal V}}

\newcommand{\Orb}{{\cal O}}
\newcommand{\R}{{\mathbb R}}

\newcommand{\td}{{\rm d}}

\makeatletter
\renewcommand{\endexample}{\hfill$\lozenge$\@endtheorem}
\makeatother

\begin{document}

\begin{center}\textbf{ {\Large Singular reduction of implicit Hamiltonian systems}\\[2ex]
Guido~Blankenstein\\
Department of Mechanical Engineering\\
Katholieke Universiteit Leuven\\
Celestijnenlaan 300 B\\
B-3001 Leuven (Heverlee)\\
Belgium\\[2ex]
Tudor~S.~Ratiu\\
Institut de Math\'ematiques Bernoulli\\
\'Ecole Polytechnique F\'ed\'erale de Lausanne\\ 
MA-Ecublens, CH-1015 Lausanne\\ 
Switzerland
}
\end{center}

\begin{abstract}
This paper develops the reduction theory of implicit Hamiltonian systems 
admitting a symmetry group at a \emph{singular} value of the momentum map. 
The results naturally extend those known for (explicit) Hamiltonian systems described by Poisson brackets.
\end{abstract}

{\bf keywords:} implicit Hamiltonian systems, Dirac structures, symmetry, reduction

\begin{center}
\begin{minipage}[t]{10cm} {\footnotesize \tableofcontents}\end{minipage}

\end{center}

\section{Introduction}
\setcounter{equation}{0}

Consider a symplectic manifold $(M,\omega)$ admitting a symmetry Lie
group $G$ (with Lie algebra denoted by $\G$) acting freely and
properly on $M$, together with a corresponding equivariant momentum
map $P: M \rightarrow \G^\ast$.  In \cite{MW74} it is shown that, at a
\emph{regular} value $\mu \in \G^\ast$ of the momentum map, the
symplectic structure on $M$ naturally reduces to a symplectic
structure $\omega_\mu$ on the reduced manifold $M_\mu =
P^{-1}(\mu)/G_\mu$, where $G_\mu$ is the coadjoint isotropy subgroup
of $G$. Furthermore, the integral curves of a Hamiltonian vector field
corresponding to a $G$-invariant Hamiltonian $H \in C^\infty(M)^G$
project to integral curves of the reduced Hamiltonian vector field
corresponding to the reduced Hamiltonian $H_\mu \in C^\infty(M_\mu)$.
This theory has been generalized in \cite{MR86} to the case of Poisson
manifolds: at a regular value $\mu$ of the momentum map, the Poisson
bracket $\{ \cdot,\cdot \}: C^\infty(M) \times C^\infty(M) \rightarrow
C^\infty(M)$ on $M$ descends to a Poisson bracket $\{ \cdot,\cdot
\}_\mu$ on the reduced phase space $M_\mu$. Again, the Hamiltonian
flow defined by a $G$-invariant Hamiltonian $H \in C^\infty(M)^G$
reduces to a Hamiltonian flow on $M_\mu$ corresponding to $H_\mu \in
C^\infty(M_\mu)$. We refer to \cite{AM78,MR99,O93} for some excellent
overviews of this theory.

Although by Sard's theorem regular reduction describes the ``generic"
case (i.e.,~the regular values of the momentum map are dense in
$\G^\ast$), interesting dynamics, such as bifurcation phenomena, occur
at the \emph{singular} values of the momentum map. Already in simple
examples like the spherical pendulum moving with angular momentum zero
(i.e.,~moving in a plane) and the Lagrange top (reducing the
gravitational $S^1$ symmetry after \emph{regular} reduction of the
internal $S^1$ symmetry corresponding to the homogeneous mass
distribution), one recognizes the need for the investigation of these
special cases. This task has been taken up in \cite{ACG91,BL97,CS91,
  CSn91, OR98, OR03, SL91} leading to the theory of so-called singular
reduction of (symplectic and Poisson) Hamiltonian systems. See
\cite{CB97} for a nice overview and some worked out examples
(including the spherical pendulum and the Lagrange top). The main
difference with the regular reduction theory is that in the case of a
singular value $\mu$ of the momentum map (i.e.,~the derivative, or
tangent map, of $P$ at points of $P^{-1}(\mu)$ is not surjective), the
reduced space $M_\mu$ is \emph{not} a manifold. Therefore symplectic
forms and (Hamiltonian) vector fields are not defined and, as a
consequence, the reduced ``Hamiltonian dynamics" cannot be written as
a system of ordinary differential equations on $M_\mu$ (as in the
regular case). However, since $M_\mu$ is still a topological space
(relative to the natural quotient topology), a reduced Poisson bracket
$\{ \cdot,\cdot\}_\mu$ on the space of (Whitney) smooth functions
$C^\infty(M_\mu)$ can still be defined.  This bracket induces a
Hamiltonian formalism that allows one to write the reduced Hamiltonian
dynamics on the singular reduced space $M_\mu$. The Hamiltonian flow
corresponding to this Hamiltonian dynamics is exactly the projection
of the regular Hamiltonian flow on $M$. Finally, in
\cite{BL97,CB97,CS91,CSn91,OR03,SL91} it has been shown that the
singular reduced space $M_\mu$ (resulting from a \emph{symplectic}
manifold $(M,\omega)$) may be stratified by symplectic manifolds,
called pieces. The stratification is by orbit type decomposition. The
Hamiltonian flow on $M_\mu$ leaves these pieces invariant and
restricts to a (regular) Hamiltonian flow on each of the pieces.

Recently, \cite{B00,BvdS01} have generalized the \emph{regular}
reduction theory for explicit Hamiltonian systems to a regular
reduction theory for \emph{implicit} Hamiltonian systems, extending
preliminary work in \cite{C90,vdS98}. Analogous to the symplectic
form or the Poisson bracket in the classical theory, the underlying
geometric structure of an implicit Hamiltonian system is that of a
\emph{Dirac structure}, defined as a maximally isotropic smooth vector
subbundle of $TM \oplus T^\ast M$. This structure allows to define a
Hamiltonian formalism generalizing the classical symplectic and
Poisson formalisms by including the description of \emph{algebraic
  equations}. That is, the ``Hamiltonian dynamics" corresponding to a
Dirac structure and a function $H \in C^\infty(M)$ consists of a set
of differential and algebraic equations (and is therefore called
implicit). Perhaps the most striking example of an implicit
Hamiltonian system is that of a Hamiltonian system defined by a
Poisson bracket on $M$, restricted to a submanifold of $M$ that is \emph{not}
the level set of a Casimir function (this example actually
motivated the definition of a Dirac structure in \cite{C90}).  Using
the notion of a Dirac structure introduced in \cite{C90,D93}, implicit
Hamiltonian systems were defined in
\cite{vdSM95,vdSM95b,vdSM97,MvdS97,BC97} and successfully employed in
the context of network modeling of energy conserving physical systems
such as mechanical systems with (non)holonomic kinematic constraints,
electrical LC circuits, electromechanical systems. We refer to
\cite{vdS99b} and references therein for more information; see also
\cite{B00} for a detailed historical account. (Note that the
Hamiltonian systems defined in \cite{C90,D93} are not truly implicit,
but rather \emph{explicit} systems.) The study of symmetries and
reduction of implicit Hamiltonian systems evolved from preliminary
results in \cite{C90,D93,vdS98} and led to a reduction theory
for this class of systems described in \cite{B00,BvdS01}. There it was
shown that a Dirac structure on $M$, admitting a symmetry Lie group
$G$ with corresponding equivariant momentum map $P$, reduces to a
Dirac structure on the reduced manifold $M_\mu$ if $\mu$ is a
\emph{regular} value of the momentum map. Furthermore, the
(projectable) integral curves of the implicit system corresponding to
a $G$-invariant function $H \in C^\infty(M)^G$ project to integral
curves of the reduced implicit Hamiltonian system defined by the
reduced Dirac structure and the reduced Hamiltonian $H_\mu$. The
theory generalizes the classical regular reduction theory for
symplectic and Poisson Hamiltonian systems, as well as the recently
developed reduction theories for constrained mechanical systems. 
Section \ref{sec:sr} briefly recalls the main results of
\cite{B00,BvdS01}; for a discussion of the connection with the
reduction theory for constrained mechanical systems we refer to
\cite{B00,BvdS01}.

The goal of this paper is to develop reduction theory for implicit
Hamiltonian systems at \emph{singular} values of the momentum map.  We
restrict our attention to a fairly general class of Dirac structures,
described by a generalized Poisson bracket and a distribution of
derivations (i.e.,~vector fields) on smooth functions.  We consider
the special subclass of symmetries that are symmetries of both the
generalized Poisson bracket and the distribution. Using these
ingredients, we prove that one can define a so-called
\emph{topological Dirac structure} on the singular reduced space
$M_\mu$ (where for easiness of exposition we will take $\mu = 0$),
representing the reduced Dirac structure. This topological Dirac
structure, whose construction implicitly uses the Sikorski
differential spaces (see \cite{CSn91, S}), defines a Hamiltonian
formalism on the reduced space $M_\mu$ and the dynamics corresponding
to an implicit Hamiltonian system with Hamiltonian $H \in
C^\infty(M_\mu)$ are described. It is shown that the (projectable)
integral curves of the implicit Hamiltonian system on $M$
corresponding to a $G$-invariant function $H \in C^\infty(M)^G$
project to ``integral curves" of the Hamiltonian dynamics defined on
the singular reduced space. The orbit type decomposition defines a
stratification of the singular reduced space. It is shown that, under
certain conditions, the topological Dirac structure restricts to
regular reduced Dirac structures on the pieces (which are always
manifolds) and that the Hamiltonian flow restricts to a regular
Hamiltonian flow on each of these pieces.

The paper is organized as follows. Section \ref{sec:ihs} gives a brief
introduction to Dirac structures and implicit Hamiltonian systems. The
basic results concerning symmetries and regular reduction of implicit
Hamiltonian systems are recalled in section \ref{sec:sr}. Section
\ref{sec:sing} describes in a purely topological way how to reduce an
implicit Hamiltonian system admitting a symmetry group. The result is
a so-called topological Dirac structure on the singular reduced space
(in general, not a manifold). It is shown that if the symmetry group
acts regularly and the value of the momentum map is regular, and hence
the reduced space is a smooth manifold, then the singular reduced
implicit Hamiltonian system equals the regular reduced implicit
Hamiltonian system as defined in section \ref{sec:sr}. Section
\ref{sec:sing_dyn} discusses the dynamics of a singular
reduced implicit Hamiltonian system. It is shown that the
``projectable" solutions of the original system project to solutions
of the singular reduced system. Section \ref{sec:orbitdecomp}
describes the decomposition of the singular reduced space into smooth
manifolds called pieces. The singular reduced implicit Hamiltonian
system restricts to regular reduced implicit Hamiltonian systems on
the pieces. Section \ref{sec:concl} contains some concluding remarks.

\section{Implicit Hamiltonian systems}\label{sec:ihs}
\setcounter{equation}{0}

Let $M$ be a smooth $n$-dimensional manifold and let $TM \oplus T^\ast
M$ denote the vector bundle whose fiber at $x \in M$ is $T_xM \times
T_x^\ast M$. (Throughout this paper all geometric objects are assumed
to be smooth, so when manifolds, vector bundles, sections are
mentioned, they are all smooth.) A Dirac structure on $M$ is defined
as follows.

\begin{definition}\label{def:ds}
  A smooth vector subbundle $\D \subset TM \oplus T^\ast M$ is called
  a \emph{Dirac structure} if for every fiber $\D(x) \subset T_xM
  \times T_x^\ast M, \; x\in M,$ one has $\D(x) = \D^\perp (x)$, where
\begin{equation}\label{def:dirac}
\D^\perp (x) = \{ (w,w^\ast) \in T_xM \times T_x^\ast M \mid \langle v^\ast,w \rangle + \langle w^\ast, v \rangle = 0, \; \forall (v,v^\ast) \in \D(x) \}.
\end{equation}
(Here $\langle \cdot,\cdot \rangle$ denotes the duality pairing
between $TM$ and $T^\ast M$.)
\end{definition}

Notice that $\D$ being a vector subbundle of $TM \oplus T^\ast M$
implies that its fibers all have the same dimension, i.e., $\dim \D(x)
= \dim \D(x^\prime), \; \forall x,x^\prime \in M$. In particular, if
$\D$ is a Dirac structure then $\dim \D(x) = n, \; \forall x \in M$.
Furthermore, $\D(x) = \D^\perp (x), \; x \in M,$ implies that
\begin{equation}\label{ortho}
\langle v^\ast,v \rangle = 0, \; \forall (v,v^\ast) \in \D(x).
\end{equation}  

\begin{remark} 
{\rm In \cite{DvdS99} a \emph{constant Dirac structure} on a vector
    space $\V$ is defined as a vector subspace $\D \subset \V \times
    \V^\ast$ such that $\D = \D^\perp$. An equivalent way of writing
    Definition \ref{def:ds} is therefore: a Dirac structure on a
    manifold $M$ is a smooth vector subbundle $\D \subset TM \oplus
    T^\ast M$ such that each fiber $\D(x), \;x \in M,$ is a constant
    Dirac structure on $T_xM$.  
} \quad $\blacklozenge$
\end{remark}

In \cite{C90} there is yet another slightly different definition of a
Dirac structure. Denote by $\mathfrak{X}_{loc}(M)$ (respectively
$\mathfrak{X}(M)$) the space of local (respectively global) smooth
sections of $TM$, that is, these are the spaces of smooth local
(respectively global) vector fields on $M$. Similarly,
$\Omega^k_{loc}(M)$ and $\Omega^k(M)$ denote the spaces of smooth
local and global $k$-forms on $M$. The spaces of smooth local and
global sections of the vector subbundle $\D \subset TM \oplus T^\ast
M$ are denoted by $\mathfrak{D}_{loc}$ and $\mathfrak{D}$
respectively. Throughout, let $X, Y \in \mathfrak{X}_{loc}(M)$ and
$\alpha, \beta \in \Omega^1_{loc}(M)$. Define a pairing on smooth
sections of $TM \oplus T^\ast M$ by
\begin{equation}\label{eq:pairing}
\langle \langle (X,\alpha),\:(Y,\beta) \rangle \rangle = \langle \alpha,Y \rangle + \langle \beta,X \rangle, \; \mbox{for}\;(X,\alpha),\:(Y,\beta) \in TM \oplus T^\ast M.
\end{equation}
According to \cite{C90}, a Dirac structure on $M$ is a smooth vector
subbundle $\D \subset TM \oplus T^\ast M$ such that
\begin{itemize}
\item[1.] $\D$ is isotropic, i.e.,~for every two (local) sections
  $(X,\alpha),\:(Y,\beta) \in \mathfrak{D}_{loc}$: $\langle \langle
  (X,\alpha),\:(Y,\beta) \rangle \rangle = 0$;
\item[2.] $\D$ is maximal, i.e.,~if $(Y,\beta)$ is a (local) section
  of $TM\oplus T^\ast M$ such that $\langle \langle
  (X,\alpha),\:(Y,\beta) \rangle \rangle = 0, \; \forall (X,\alpha)
  \in \mathfrak{D}_{loc}$, then $(Y,\beta) \in \mathfrak{D}_{loc}$.
\end{itemize}

It is easily shown that this definition given in
\cite{C90} and Definition \ref{def:ds} are equivalent. Indeed,
since $\D$ is a smooth vector subbundle, every
$(v,v^\ast) \in \D(x)$ can be extended to a local section $(X,\alpha)
\in \mathfrak{D}_{loc}$. Furthermore, $\D$ being a smooth vector
subbundle implies that also $\D^\perp \subset TM \oplus T^\ast M$
(with fibers $\D^\perp(x)$) is a smooth vector subbundle and therefore
also every $(w,w^\ast) \in \D^\perp (x)$ can be extended to a local
section $(Y,\beta)$ of $\D^\perp$. Elementary linear algebra shows
that $\D$ is isotropic if and only if $\D \subset \D^\perp$ and that
maximal isotropy is equivalent to equality in this inclusion or to the
fact that $\operatorname{dim} \D(x) = n$ for all $x\in M$ (see,
e.g., \cite {AM78}, Section 5.3 for this type of argument).
  
\medskip

Corresponding to a Dirac structure a number of (co-)distributions may
be identified. Recall that a \textit{distribution\/} $\Delta$ on a
manifold $M$ is a map which assigns to each $x \in M$ a linear
subspace $\Delta(x)$ of the tangent space $T_xM$. $\Delta$ is called a
\textit{smooth distribution\/} if around any point these subspaces are
spanned by a set of smooth vector fields, i.e., there exist
$X_1,\dots,X_k \in \mathfrak{X}_{loc}(M)$ such that $\Delta(x) =
\mbox{span}\; \{ X_1(x),\dots,X_k(x) \}$. The distribution $\Delta$ is
called \textit{constant dimensional\/} if the dimension of the linear
subspace $\Delta(x) \subset T_xM$ does not depend on the point $x \in
M$. Notice that if $\Delta$ is a smooth constant dimensional
distribution on $M$, then it defines a smooth vector subbundle (also
denoted by $\Delta$) of the tangent bundle $TM$, with fibers
$\Delta(x)$, $x \in M$. Analogously, a \textit{codistribution\/}
$\Gamma$ is defined as a map which assigns to each $x \in M$ a
linear subspace $\Gamma(x)$ of the cotangent space $T_x^\ast M$.
Smoothness and constant dimensionality are defined in the same way as
for distributions. A smooth constant dimensional codistribution
defines a smooth vector subbundle of the cotangent bundle $T^\ast M$.

Now, any Dirac structure $\D$ naturally defines a distribution
$\Delta$, with fibers given by\footnote{In the literature on implicit
  Hamiltonian systems this distribution is usually denoted by
  $\mathsf{G}_0$ and sometimes called the characteristic distribution.
  However, in order to avoid confusion with notation defined later on
  in this paper we decided to adopt a different notation here.}
\begin{equation}
\Delta(x):=\{X(x) \mid X \in \mathfrak{X}_{loc}(M), (X,0) \in \mathfrak{D}_{loc}\},
\end{equation}
and a codistribution $\Gamma$, whose fibers are defined
by\footnote{This codistribution is usually denoted by $\mathsf{P}_1$
  in the literature.}
\begin{equation}\label{gamma}
\Gamma(x) := \{ \alpha(x) \mid  \alpha \in \Omega_{loc}^1 (M), \exists X \in \mathfrak{X}_{loc}(M), (X,\alpha) \in \mathfrak{D}_{loc} \}.
\end{equation}
Since $\D$ is isotropic it follows that $\Delta(x) \subset
\Gamma^\circ(x)$, where $\Gamma^\circ(x)$ denotes the annihilating
vector subspace of $\Gamma(x)$ in $T_xM$, that is, $\Gamma^\circ(x) =
\{ v \in T_xM \mid \langle v^\ast,v \rangle = 0, \;\forall v^\ast \in
\Gamma(x) \}$. Equivalently, $\Gamma(x) \subset \Delta^\circ(x)$,
where $\Delta^\circ (x)$ denotes the annihilating vector subspace of
$\Delta(x)$ in $T_x^\ast M$, that is, $\Delta^\circ (x): = \{ v^\ast
\in T^\ast _x M \mid \langle v^\ast, v \rangle = 0, \;\forall v \in
\Delta(x) \}$. Furthermore, if $\Gamma$ is constant dimensional, and
hence defines a vector subbundle of $T^\ast M$, it follows by maximal
isotropy of $\D$ that $\Delta(x) = \Gamma^\circ(x)$ (equivalently,
$\Gamma(x) = \Delta^\circ(x)$). Notice that in this case also $\Delta$
is constant dimensional and hence defines a vector subbundle of $TM$.
\medskip

A special property of a Dirac structure is defined in the following
\begin{definition}\label{def:closed}
  A Dirac structure $\D$ is called \emph{closed}, or
  \emph{integrable}, if for all (local) sections $(X_1,\alpha_1)$,
  $(X_2,\alpha_2), \:(X_3,\alpha_3) \in \mathfrak{D}_{loc}$
\begin{equation}
\langle L_{X_1} \alpha_2, X_3 \rangle + \langle L_{X_2} \alpha_3, X_1 \rangle + \langle L_{X_3} \alpha_1, X_2 \rangle = 0.
\end{equation}

Equivalently \cite{C90,D93,DvdS99}, $\D$ is closed if and only if for
all $(X_1,\alpha_1), \:(X_2,\alpha_2) \in \mathfrak{D}_{loc}$
\begin{equation}\label{eq:closed}
\left( [X_1,X_2], L_{X_1} \alpha_2 - L_{X_2} \alpha_1 + \td \langle \alpha_1,X_2 \rangle \right) \in \mathfrak{D}_{loc}.
\end{equation}
\end{definition}

The notation $L_X$ is reserved for the Lie derivative operator (acting
on any type of tensor field) defined by the (local) vector field $X$
on $M$.

It is easy to see that the graph of a symplectic form $\omega: TM
\rightarrow T^\ast M$ or the graph of the skew-symmetric vector bundle
map $J: T^\ast M \rightarrow TM$ induced by a Poisson bracket on $M$
defines a Dirac structure on $M$. As customary, we will call both the
bundle map $J$ and the two-tensor defined by $\{\cdot, \cdot \}$ the
Poisson structure on $M$. Closedness of $\D$ corresponds in these two
cases to the condition that $\omega$ is a closed two-form,
respectively, the Poisson bracket satisfies the Jacobi identity.
\medskip

In this paper we will concentrate on a rather frequently occurring type
of Dirac structure defined as follows. Let $\{ \cdot,\cdot \}:
C^\infty(M) \times C^\infty(M) \rightarrow C^\infty(M)$ be a
\emph{generalized Poisson bracket} on $M$. That is, $\{ \cdot,\cdot
\}$ is skew-symmetric, bilinear, and satisfies the Leibniz property.
Denote the corresponding vector bundle map by $J: T^\ast M \rightarrow
TM$ (i.e., $J(dH,dF) = \{H,F\}, \; H,F \in C^\infty(M)$); recall that
$J$ is skew-symmetric. Note that we do not require $\{ \cdot,\cdot
\}$ to satisfy the Jacobi identity and neither that $J$ have constant
rank. Moreover, given a subbundle $\Delta$ of $TM$ (i.e.,
a smooth constant dimensional distribution $\Delta$ on $M$), it is
easy to see that the vector subbundle $\D \subset TM \oplus T^\ast
M$ with fiber
\begin{equation}\label{ds1}
\D(x) = \{ (v,v^\ast) \in T_xM \times T_x^\ast M \mid v - J(x) v^\ast \in \Delta(x), \; v^\ast \in \Delta^\circ (x) \}
\end{equation}
defines a Dirac structure on $M$. In terms of its local sections this is expressed as
\begin{equation}\label{ds2}
\mathfrak{D}_{loc} = \{ (X,\alpha) \in \mathfrak{X}_{loc}(M) \oplus \Omega^1_{loc}(M)  \mid X - J\alpha \;\mbox{is a local section of} \; \Delta, \; \alpha \;\mbox{ is local section of}\; \Delta^\circ \},
\end{equation}
where $\Delta^\circ$ denotes the vector subbundle of $T^\ast M$ whose fiber at $x \in M$ equals $\Delta^\circ (x)$. 

\begin{remark}
  {\rm In \cite{DvdS99} it is shown that, under a mild
    constant dimensionality assumption, \emph{every} Dirac structure
    can be written in the form (\ref{ds1}) or, equivalently,
    (\ref{ds2}). Indeed, if $\D$ is an arbitrary Dirac structure on
    $M$, define the codistribution $\Gamma$ as in (\ref{gamma}), and
    assume that $\Gamma$ is constant dimensional and hence defines a
    vector subbundle of $T^\ast M$. Then there exists a well defined
    (see \cite{DvdS99}) skew-symmetric vector bundle map $J(x):
    \Gamma(x) \subset T_x^\ast M \rightarrow (\Gamma(x))^\ast \subset
    T_xM, \; x \in M$, defined by
\begin{equation}
J(x) v^\ast = \bar{v} \in (\Gamma(x))^\ast, \; v^\ast \in \Gamma(x),
\end{equation}
where $\bar{v} \in (\Gamma(x))^\ast$ is the restriction of some $v \in
T_xM$ to $\Gamma(x) \subset T_x^\ast M$ which satisfies the condition
that $(v,v^\ast ) \in \D(x)$. Notice that the kernel of $J(x)$ is
given by the codistribution $\Gamma_0$ with fibers defined
by\footnote{This codistribution is usually denoted by $\mathsf{P}_0$
  in the literature.}
\begin{equation}
\Gamma_0 (x) := \{ \alpha(x) \mid \alpha \in \Omega_{loc}^1 (M), (0,\alpha) \in \mathfrak{D}_{loc} \}.
\end{equation}
Then the Dirac structure $\D$ takes the form (\ref{ds1}), or,
equivalently, (\ref{ds2}), with $\Delta = \Gamma^\circ$.

The map $J$ may be locally extended to a generalized Poisson structure
$J: T^\ast M \rightarrow TM$. 
} \quad $\blacklozenge$ 
\end{remark}  
 
Notice that, in general, the Dirac structure defined in (\ref{ds1}) is
not closed. Although closedness does not play an important role in the
rest of the paper we remark for completeness that (\ref{ds2}) defines
a closed Dirac structure if and only if (see \cite{DvdS99})
\begin{itemize}
\item[1.] $\Delta$ is involutive,
\item[2.] the bracket $\{ \cdot,\cdot \}$ restricted to the set of
  \emph{admissible functions} $\A_D := \{ H \in C^\infty(M) \mid dH \in
  \Delta^\circ \}$ defines a Poisson bracket on $\A_D$ (that is, the
  Jacobi identity holds).
\end{itemize}

\begin{remark}\label{rem:algebroids}{\rm
    Before leaving our introduction to Dirac structures and moving on
    to the description of implicit Hamiltonian systems, we would like
    to mention the following generalization. In \cite{LWX98} a Dirac
    structure is defined as a maximal isotropic subbundle $\D \subset
    \A \oplus \A^\ast$, where the pair $(\A,\A^\ast)$ defines a
    \emph{Lie bialgebroid} over a smooth manifold $M$. The isotropy is
    defined with respect to the natural pairing $\langle \langle ,
    \rangle \rangle$ defined analogously as in (\ref{eq:pairing}) by
    $\A$ and its dual $\A^\ast$. If we take the special case $\A=TM$
    and dually $\A^\ast = T^\ast M$, then we are back at the
    definition given earlier in this paper. For more information on
    this generalization we refer to \cite{LWX98} and the references
    therein. We remark that \cite{LWX98} require the Dirac structure
    to be closed. In their terminology, \emph{Dirac structures} are
    always closed, while Dirac structures \emph{not} satisfying
    condition (\ref{eq:closed}) are called \emph{almost Dirac
      structures}. In this paper however, we prefer to use the
    terminology as introduced above, calling a maximally isotropic
    subbundle of $TM \oplus T^\ast M$ a Dirac structure and adding the
    prefix \emph{closed} if (and only if) the conditions of Definition
    \ref{def:closed} are satisfied.  }\hfill $\blacklozenge$
\end{remark}

Now we turn to the definition of an implicit Hamiltonian system.
Consider a Dirac structure $\D$ on $M$ and a smooth function $H \in
C^\infty(M)$, called the Hamiltonian or energy function. Then the
three-tuple $(M,\D,H)$ defines an implicit Hamiltonian system as
follows:
\begin{definition}\label{def:ihs}
  The \emph{implicit Hamiltonian system} $(M,\D,H)$ is defined as a
  set of smooth time functions $\{ x(t) \mid x:\R \rightarrow M
  \text{~of~class~} C^\infty\}$ (called \emph{solutions}) satisfying
  the condition
\begin{equation}\label{def:ihs-sol}
(\dot{x}(t), dH(x(t))) \in \D(x(t)), \; \forall t.
\end{equation}
\end{definition}

Equations (\ref{ortho}) and (\ref{def:ihs-sol}) imply that implicit Hamiltonian systems are \emph{energy conser\-ving}, i.e.,
\begin{equation}\label{ec}
\frac{dH}{dt}(x(t)) = \langle dH(x(t)),\dot{x}(t) \rangle = 0, \; \forall t.
\end{equation}

If $\D$ is the graph of a symplectic form $\omega$ or of a Poisson
structure $J: T^\ast M \rightarrow TM$ then the above definition
becomes that of a classical (explicit) symplectic or Poisson
Hamiltonian system. On the other hand, if $\D$ is defined by
(\ref{ds1}) then the system includes the \emph{algebraic constraints}
\begin{equation}
dH(x(t)) \in \Delta^\circ (x(t)), \; \forall t.
\end{equation}
Thus all solutions of the implicit Hamiltonian system
necessarily lie in the \emph{constraint manifold}
\begin{equation}
M_c = \{ x \in M \mid dH(x) \in \Delta^\circ (x) \}.
\end{equation}

Since the implicit Hamiltonian system defines a set of differential
and algebraic equations, there is not an existence and
uniqueness result as one has for solutions of classical Hamiltonian
systems described by ordinary differential equations. In particular,
not every point $x_0 \in M_c$ necessarily lies on the trajectory of
some solution $x(t)$ of the system, and neither are the solutions
through a point $x_0 \in M_c$ (if existing) necessarily unique (this
happens, for instance, if the Lagrange multipliers corresponding to the
algebraic constraints cannot be solved uniquely). In the sequel we
will not investigate these problems. Instead we will study the
reduction of the underlying Dirac structure in the presence of
symmetries (defined later on) and show that certain ``projectable"
solutions (if existing) will project to solutions of an implicit
Hamiltonian system on the reduced space.

The problem of existence and uniqueness of solutions to implicit
systems is an important and active area of research and will not be
touched upon here. We would only like to mention the special case of
so-called \emph{index 1} systems. Consider the implicit Hamiltonian
system defined by the Dirac structure (\ref{ds1}) and the Hamiltonian
function $H \in C^\infty(M)$. Let the vector subbundle $\Delta$ be
(locally) written as the span of the independent vector fields
$g_1,\dots,g_m$. Then the constraint manifold can be written as
\begin{equation}
M_c = \{ x \in M \mid L_{g_j} H(x)=0, \; j=1,\dots,m \}.
\end{equation}
Now \emph{assume} that the constraints are of index 1, that is, the matrix
\begin{equation}
[L_{g_i}L_{g_j} H(x)]_{i,j=1,\dots,m}
\end{equation}
is nonsingular for all $x \in M_c$. In that case one can restrict the
implicit Hamiltonian system $(M,\D,H)$ to an \emph{explicit}
Hamiltonian system on $M_c$, defined by a (possibly nonintegrable)
Poisson bracket on $M_c$, see \cite{vdSM95b,B00}. Its corresponding
dynamics is thus given by a set of ordinary differential
equations on $M_c$. Standard existence and uniqueness results now
yield that through every point $x_0 \in M_c$ there is a unique
solution of the implicit Hamiltonian system (restricted to $M_c$).

\begin{example}\label{ex:ms}
  {\rm As an important example of implicit Hamiltonian systems we like
    to mention the class of mechanical systems with kinematic
    constraints. These systems are described by implicit Hamiltonian
    systems $(M,\D,H)$ with $\D$ of the form (\ref{ds1}). Here, $H$ is
    the total energy of the system, (the sum of kinetic and potential
    energies for classical mechanical systems) and the phase space $M
    = T^\ast Q$ is the cotangent bundle of the configuration space
    $Q$; local coordinates are denoted, as usual, by $(q,p) \in T^\ast
    Q$. The Poisson bracket $\{ \cdot,\cdot \}$ is the standard
    Poisson bracket corresponding to the canonical symplectic form
    $\omega = dp \wedge dq$ on $T^\ast Q$ (i.e.,~the associated
    Poisson structure $J$ is given by $\omega^{-1}$). Finally, if we
    assume the kinematic constraints to be linear in the velocities
    (e.g.,~non-slipping constraints of a rolling wheel) then they can
    be (locally) written in the form
\begin{equation}\label{ex:cons}
A^T(q)\dot{q} = A^T(q) \frac{\partial H}{\partial p}(q,p) = 0.
\end{equation}
By d'Alembert's principle, the constraints (\ref{ex:cons}) generate
constraint forces of the form $F_{c} = A(q)\lambda$, where $\lambda$
are the Lagrange multipliers. It follows that the distribution
$\Delta$ is (locally) described by the image of the matrix
\begin{equation}\label{mechDelta}
\left[ \begin{array}{c} 0\\ A(q) \end{array} \right].
\end{equation}

If the kinetic energy is defined by a positive definite metric on $Q$,
then the constraints are of index 1, i.e.,~the Lagrange multipliers
$\lambda$ can be solved uniquely. In this case the constrained
mechanical system on $T^\ast Q$ can be written as an
\emph{unconstrained} generalized Hamiltonian system on $M_c$. In
\cite{vdSM94} it is shown that the corresponding Poisson bracket on
$M_c$ satisfies the Jacobi identity if and only if the kinematic
constraints are \emph{holonomic}.  
}
\end{example}

\section{Symmetries and regular reduction}\label{sec:sr}
\setcounter{equation}{0}

In this section we recall some of the results in \cite{B00,BvdS01}
concerning symmetries and reduction of implicit Hamiltonian systems.
These results act as a reference for the results obtained in this
paper and will be specifically needed in section \ref{sec:orbitdecomp}
to show that the singular reduced Dirac structure restricts to regular
reduced Dirac structures on the pieces corresponding to the orbit type
decomposition of the singular reduced space $M_0$. We refer to
\cite{B00,BvdS01} for a detailed treatment of symmetries and reduction
of implicit Hamiltonian systems. We stress that, unless specifically
stated otherwise, the results in this section are valid for arbitrary
Dirac structures, not necessarily of the form (\ref{ds1}).

\begin{definition}\label{def:sym}
  A smooth vector field $Y$ on $M$ is called a \emph{symmetry} of the
  Dirac structure $\D$ if for every (local) section $(X,\alpha)$ of $
  \D$, one also has that $(L_YX,L_Y\alpha)$ is a (local) section of
  $\D$.  $Y$ is called a \emph{symmetry} of the implicit Hamiltonian
  system $(M,\D,H)$ if $Y$ is a symmetry of $\D$ and a symmetry of
  $H$, i.e.,~$L_Y H=0$.\footnote{This is called a \emph{strong}
    symmetry in \cite{B00,BvdS01}.}
\end{definition}

The above notion of symmetry generalizes the classical notion of
symmetry in symplectic or Poisson Hamiltonian systems. Indeed, if $\D$
is the graph of a symplectic form $\omega$ then $Y$ is a symmetry of
$\D$ if and only if $L_Y \omega=0$. Likewise, if $\D$ is the graph of
the skew-symmetric vector bundle map $J$ defining a Poisson structure
on $M$, then $Y$ is a symmetry of $\D$ if and only if the integral
flow of $Y$ is a Poisson map, i.e.,~$Y$ is a derivation of the Poisson
bracket: $L_Y\{F_1, F_2\} = \{L_YF_1, F_2\} + \{F_1 , L_YF_2\}$ for
any $F_1, F_2 \in C^\infty (M)$.

In the next sections we will consider a special subclass of symmetries
defined by the following
\begin{proposition}\label{prop:sym}
  Consider a Dirac structure $\D$ of the type defined in (\ref{ds1}).
  If the vector field $Y$ on $M$ is a derivation of the generalized
  Poisson bracket (equivalently, $L_Y J=0$) and $L_Y Z$ is a local
  section of $\Delta$ whenever $Z$ is a local section of $\Delta$,
  then $Y$ is a symmetry of $\D$.
\end{proposition}

In particular this means that we restrict our attention to the case
where $Y$ is a symmetry of the generalized Poisson bracket as well as
a symmetry of the vector subbundle $\Delta \subset TM$. These kinds of
symmetries often arise in constrained mechanical systems, see also
Remark \ref{rem:mechsym} later on.

Specifically we will consider \emph{ Lie algebra symmetries}, defined
as follows. Recall that a (left) Lie algebra action on a manifold $M$
is a Lie algebra anti-homomorphism $\xi \in \mathfrak{g} \mapsto \xi_M
\in \mathfrak{X}(M)$ such that the map $(x, \xi) \in M \times
\mathfrak{g} \mapsto \xi_M(x) \in T_xM$ is smooth. Then $\mathfrak{g}$
is a \textit{symmetry Lie algebra\/} of $\D$ if $\xi_M$ is a symmetry
of $\D$ for every $\xi \in \mathfrak{g}$. In particular, if the Dirac
structure is of the type defined in (\ref{ds1}), the criterion in
Proposition \ref{prop:sym} applies. Similarly, $\mathfrak{g}$ is a
\textit{symmetry Lie algebra\/} of an implicit Hamiltonian system
$(M,\D,H)$ if each $\xi_M$ is a symmetry of $\D$ and also a symmetry
of $H$, i.e.,~$L_{\xi_M}H=0$.

Lie algebra symmetries are often induced by Lie group actions. If $G$
is a Lie group with Lie algebra $\mathfrak{g}$ and $\phi: G \times M
\rightarrow M$ is a smooth left action of $G$ on the manifold $M$, the
\textit{infinitesimal generator\/} of the action associated to $\xi\in
\mathfrak{g}$ defined by
\begin{equation}
\xi_M (x) = \frac{d}{dt} \phi(\exp(\xi t),x) \vert_{t=0} \in T_xM, \; x \in M,
\end{equation}
induces a left Lie algebra action of $\mathfrak{g}$ on $M$. Then $G$
is said to be a \textit{symmetry Lie group\/} of $\D$ if
$\mathfrak{g}$ is a symmetry Lie algebra of $\D$. Similarly, $G$ is a
\textit{symmetry Lie group\/} of the implicit Hamiltonian system $(M,
\D, H)$ if $\mathfrak{g}$ is a symmetry Lie algebra of this implicit
Hamiltonian system.

We turn now to the analysis of the regular reduction process of Dirac
structures and implicit Hamiltonian systems. We start by explaining
how an implicit Hamiltonian system on $M$ can be restricted to an
implicit Hamiltonian system on a submanifold $N$ of $M$. Let $\D$ be a
Dirac structure on $M$ and let $N \subset M$ be a submanifold of $M$.
Following \cite{C90} define for each $x \in N$ the map $\sigma (x):
T_xN \times T_x^\ast M \rightarrow T_xN \times T_x^\ast N, \; x \in
N,$ by $\sigma (x) (v,v^\ast) = (v,v^\ast\vert_{T_xN})$, where
$v^\ast\vert_{T_xN}$ means the restriction of the covector $v^\ast \in
T_x^\ast M$ to the subspace $T_xN \subset T_xM$. Define a vector
subspace of $T_xN \times T_x^\ast N$ by
\begin{equation}\label{DN}
\D_N(x)= \sigma(x) \left( \D(x) \cap (T_xN \times T_x^\ast M) \right), \; x \in N.
\end{equation}
It is clear that $\D_N(x) \subset \D_N^\perp(x), \; x \in N$. To prove
the reverse inclusion, suppose that $(w,w^\ast) \in \D_N^\perp (x)
\subset T_xN \times T_x^\ast N$, i.e., $\langle v^\ast,w \rangle +
\langle w^\ast, v \rangle = 0, \; \forall (v,v^\ast) \in \D_N(x)$.
Then $(v, v^\ast) = \sigma(x)(v, \bar{v}^\ast)$, where $\bar{v}^\ast
\in T_x^\ast M$, $(v, \bar{v}^\ast) \in \D(x)$,
$\bar{v}^\ast\vert_{T_xN} = v^\ast$ and since $v, w \in T_xN$, one
gets
\[
0=\langle v^\ast,w \rangle + \langle w^\ast, v \rangle = \langle \bar{v}^\ast,w \rangle + \langle \bar{w}^\ast, v \rangle,
\]
where $\bar{w}^\ast$ is an arbitrary extension of $w^\ast$ to $T_xM$.
Since this relation holds for all $(v, \bar{v}^\ast) \in \D(x)$ with
$v \in T_xN$, this implies that
\begin{align}\nonumber
(w,\bar{w}^\ast) &\in [ \D(x) \cap (T_xN \times T_x^\ast M) ]^\perp  = \D^\perp (x) + (T_xN \times T_x^\ast M)^\perp\\
&= \D(x) + ( \{0\} \times T_xN^\circ),
\end{align}
so there exists a $\bar{u}^\ast \in T_xN^\circ \subset T_x^\ast M$
such that $(w,\bar{w}^\ast + \bar{u}^\ast ) \in \D(x)$. However, since
$w \in T_xN$ and $\sigma(x)(w, \bar{w}^\ast + \bar{u}^\ast) = (w,
(\bar{w}^\ast + \bar{u}^\ast)\vert_{T_xN}) = (w, w^\ast)$, it follows
that $(w,w^\ast) \in \D_N(x)$, which shows that $\D_N^\perp(x) \subset
\D_N(x)$. Assuming that the dimension of $\D(x) \cap (T_xN \times
T_x^\ast M)$ is independent of $x \in N$, that is, that $\D \cap (TN
\times T^\ast M)$ is a vector subbundle of $TN \times T^\ast M$, it
follows that $\sigma$ is a vector bundle map and hence that $\D_N$ is
a vector subbundle of $TN \times T^\ast N$.  So we have proved (a slightly
rewritten version of \cite{C90}):

\begin{proposition}\label{prop:DN}
  Consider a Dirac structure $\D$ on $M$ and let $N$ be a submanifold
  of $M$. Assume that $\D(x) \cap (T_xN \times T_x^\ast M), \; x \in
  N,$ has constant dimension on $N$. Then the bundle $\D_N$ with
  fibers defined by (\ref{DN}) is a Dirac structure on $N$. This is
  called the \emph{restriction} of $\D$ to $N$.
\end{proposition}

In order to do computations it is convenient to describe the
restricted Dirac structure $\D_N$ in terms of its local sections. This
gives the following proposition (an improved version of
\cite{B00,BvdS01}). Let $\iota : N \hookrightarrow M$ denote the
inclusion map.
\begin{proposition}
  Consider a Dirac structure $\D$ on $M$ and let $N$ be a submanifold
  of $M$. Assume that $\D(x) \cap (T_xN \times T_x^\ast M), \; x \in
  N,$ has constant dimension on $N$ and let $\D_N$ denote the
  restriction of $\D$ to $N$. Then $(\bar{X},\bar{\alpha})$ is a local
  section of $\D_N$ if and only if there exists a local section
  $(X,\alpha)$ of $\D$ such that $\bar{X} \sim_\iota X$ and
  $\bar{\alpha} = \iota^\ast \alpha$. Otherwise stated, in terms of
  its local sections
\begin{equation}\label{DNsec}
(\mathfrak{D}_N)_{loc} = \{ (\bar{X},\bar{\alpha}) \in \mathfrak{X}_{loc}(N) \oplus \Omega^1_{loc} (N) \mid \exists (X,\alpha) \in \mathfrak{D}_{loc} \;\mbox{such that}\; \bar{X} \sim_\iota X \; \mbox{and}\; \bar{\alpha} = \iota^\ast \alpha \}.
\end{equation}
\end{proposition}

Here $\bar{X} \sim_\iota X$ denotes the fact that $\bar{X}$ and $X$
are $\iota$-related and $\iota^\ast$ denotes the pullback by $\iota$.
It can be shown that if $\D$ is closed, then also $\D_N$ is closed.

Now, let $(M,\D,H)$ be an implicit Hamiltonian system on $M$ and let
$N$ be a submanifold of $M$ such that the constant dimension condition
of Proposition \ref{prop:DN} is satisfied and assume that (the flow
corresponding to) the solutions of $(M,\D,H)$ leave the submanifold
$N$ invariant. Restrict the Hamiltonian $H$ to a smooth function $H_N$
on $N$, i.e., $H_N = H \circ \iota$, and define the implicit
Hamiltonian system $(N,\D_N,H_N)$ on $N$. Then:

\begin{proposition}
  Every solution $x(t)$ of $(M,\D,H)$ which is contained in $N$ is a
  solution of\linebreak[4] $(N,\D_N,H_N)$.
\end{proposition}

We remark that, in general, there is not a one-to-one correspondence
between the solutions generated by the original system $(M,\D,H)$ and
those generated by the restricted system $(N,\D_N,H_N)$. (Indeed, compare this
with the case of restriction of a symplectic form $\omega$ on $M$ to
an arbitrary submanifold $N$, leading to a nontrivial kernel for
$\omega_N$.) In case $N$ happens to be the level set of a
\emph{Casimir} function of the Dirac structure $\D$ then there
\emph{is} a one-to-one correspondence between the solutions of the
original system and of the restricted system; see \cite{B00} for more
information.

Next, we explain how an implicit Hamiltonian system on $M$ admitting a
symmetry Lie group $G$ can be projected to an implicit Hamiltonian
system on the orbit space $M/G$. Consider a Dirac structure $\D$ on
$M$ and let $G$ be a symmetry Lie group of $\D$, acting \emph{regularly} on
$M$, that is, the orbit space $M/G$ is a smooth manifold and the
canonical projection map $\pi : M \rightarrow M/G$ is a surjective
submersion.  Let $V = \ker T\pi$ denote the vertical subbundle of $TM$,
with fiber $V(x) = \mbox{span}\:\{ \xi_M(x) \mid \xi \in
\mathfrak{g}\}$ for every $x \in M$. We assume that $V+\Delta$ is a
smooth vector subbundle of $TM$, i.e., its fibers all have the same
dimension. Furthermore, define the smooth vector subbundle $E \subset
TM \oplus T^\ast M$ in terms of its local sections by
\begin{equation}
\Gamma_{loc}(E)= \{ (X,\alpha) \in \mathfrak{X}_{loc}(M) \oplus \Omega^1_{loc}(M) \mid \alpha = \pi^\ast \hat{\alpha} \;\mbox{for some}\; \hat{\alpha} \in \Omega^1_{loc} (M/G) \}
\end{equation}
(where $\Gamma_{loc}(E)$ is the space of local sections of the
subbundle $E$) and assume that $\D \cap E$ is a smooth vector
subbundle of $TM \oplus T^\ast M$, i.e., its fibers all have the same
dimension. Then we have:

\begin{proposition}\label{prop:proj}
  \cite{vdS98,BvdS01} Consider a Dirac structure $\D$ on $M$ admitting
  a symmetry Lie group $G$ acting regularly on $M$. 
 Assume that $V+\Delta$ is a smooth
  vector subbundle of $TM$ and that $\D \cap E$ is a smooth vector
  subbundle of $TM \oplus T^\ast M$. Then $\D$ projects to a Dirac
  structure $\hat{\D}$ on $\hat{M}:=M/G$, described in terms of its
  local sections by
\begin{equation}\label{Dproj}
\hat{\mathfrak{D}}_{loc} = \{(\hat{X},\hat{\alpha}) \in \mathfrak{X}_{loc} (\hat{M}) \times \Omega_{loc}^1 (\hat{M}) \mid \exists (X,\alpha) \in \mathfrak{D}_{loc} \;\mbox{such that}\; X\sim_\pi \hat{X} \;\mbox{and}\; \alpha= \pi^\ast \hat{\alpha} \}.
\end{equation}
\end{proposition}

$\hat{\D}$ is called the \emph{projection} of $\D$ to $M/G$. Again, it
can be shown that closedness of $\D$ implies closedness of $\hat{\D}$.
Let $(M,\D,H)$ be an implicit Hamiltonian system admitting a symmetry
Lie group acting regularly on $M$ such that the conditions in
Proposition \ref{prop:proj} are satisfied. The $G$-invariant function
$H$ defines a function $\hat{H} \in C^\infty(M/G)$ by $H = \hat{H}
\circ \pi$. Consider the implicit Hamiltonian system
$(M/G,\hat{\D},\hat{H})$. A $G$-\emph{projectable} solution $x(t)$ of
$(M,\D,H)$ is defined as a solution $x(t)$ of $(M,\D,H)$ for which
there exists a projectable vector field $X \in \mathfrak{X}_{loc}(M)$
(i.e.,~$X \sim_\pi \hat{X}$ for some $\hat{X} \in
\mathfrak{X}_{loc}(M/G)$) such that $\dot{x}(t) = X(x(t))$. The
following proposition was obtained in \cite{B00,BvdS01}.

\begin{proposition}
  If $x(t)$ is a $G$-projectable solution of $(M,\D,H)$ then
  $\hat{x}(t) := \pi(x(t))$ is a solution of $(M/G,\hat{\D},\hat{H})$.
  Conversely, every solution $\hat{x}(t)$ of $(M/G,\hat{\D},\hat{H})$
  is locally the projection under $\pi$ of a $G$-projectable solution
  $x(t)$ of $(M,\D,H)$.
\end{proposition}

In \cite{B00} a simple example is given showing that not every
solution of an implicit Hamiltonian system $(M,\D,H)$ admitting a
symmetry Lie group $G$ is necessarily $G$-projectable. However, if the
constraints are of index 1, and therefore the implicit Hamiltonian
system can be restricted to an explicit Hamiltonian system on the
constraint manifold $M_c$, then it can be shown that \emph{every}
solution is $G$-projectable and therefore projects to a solution of
the reduced implicit Hamiltonian system on $M/G$.

Finally, let us briefly recall some results in \cite{B00,BvdS01}
concerning the reduction of implicit Hamiltonian systems admitting a
symmetry Lie group having an associated momentum map. Consider an
implicit Hamiltonian system $(M,\D,H)$ admitting a symmetry Lie group
$G$ with Lie algebra $\mathfrak{g}$. Assume that there exists an
$Ad^\ast$-equivariant map $P: M \rightarrow \mathfrak{g}^\ast$
($Ad^\ast$ denoting the coadjoint action), called \emph{momentum map},
such that
\begin{equation}\label{momentum}
(\xi_M,dP_\xi) \in \mathfrak{D}, \; \forall \xi \in \mathfrak{g},
\end{equation}
where $P_\xi \in C^\infty(M)$ is defined by $P_\xi (x) = \langle P(x),
\xi \rangle, \; x \in M$. Notice that if $\D$ is the graph of a
symplectic form $\omega$ or a Poisson structure $J$, then this
corresponds to the classical definition of a momentum map. Assuming
that $\mu \in \mathfrak{g}^\ast$ is a \emph{regular} value of $P$, it
follows that the level set $P^{-1}(\mu)$ is a closed submanifold of
$M$. Since the Hamiltonian is $G$-invariant, the solutions of
$(M,\D,H)$ leave the level set $P^{-1}(\mu)$ invariant. Thus by
(\ref{def:ihs-sol}), (\ref{momentum}), and the identity
$\D(x)=\D^\perp(x)$ it follows that
\begin{equation}\label{integral}
\frac{dP_\xi}{dt}(x(t)) = \langle dP_\xi(x(t)), \dot{x}(t) \rangle = - \langle dH,\xi_M \rangle (x(t)) = 0, \; \forall t, \; \forall \xi \in \mathfrak{g}.
\end{equation}
In other words, $P$ is a \emph{first integral} of the implicit
Hamiltonian system $(M,\D,H)$. Assuming that the conditions in Proposition
\ref{prop:DN} hold, we can restrict the implicit Hamiltonian
system $(M,\D,H)$ to an implicit Hamiltonian system $(N,\D_N,H_N)$ on
$N=P^{-1}(\mu)$. The system $(N,\D_N,H_N)$ admits the symmetry Lie
group $G_\mu := \{g\in G \mid Ad^\ast_g \mu = \mu \}$. Assume that
$G_\mu$ acts \emph{regularly} on $N$, that is, $N/G_\mu$ is a smooth
manifold with the canonical projection a surjective submersion, and
assume that the conditions in Proposition \ref{prop:proj} are
satisfied. Then we can project the implicit Hamiltonian system
$(N,\D_N,H_N)$ to an implicit Hamiltonian system
$(M_\mu,\D_\mu,H_\mu)$, where $M_\mu = N/G_\mu = P^{-1}(\mu)/G_\mu$ is
the regular reduced space and $H_\mu \in C^\infty(M_\mu)$ defined by
$H_\mu \circ \pi = H_N$ is the reduced Hamiltonian. The implicit
Hamiltonian system $(M_\mu,\D_\mu,H_\mu)$ is called the \emph{reduced}
implicit Hamiltonian system corresponding to $(M,\D,H)$ and $\D_\mu$
is called the reduced Dirac structure. If $\D$ is the graph of a
symplectic form $\omega$, then $\D_\mu$ is precisely the graph of the
Marsden-Weinstein reduced symplectic form $\omega_\mu$. Likewise, if
$\D$ is the graph of a Poisson structure $J$ on $M$, then $\D_\mu$ is
the graph of the reduced Poisson structure $J_\mu$ \cite{MR86}.
Notice, however, that contrary to the above mentioned classical
reduction results, closedness of the Dirac structure is \emph{not}
required for the reduction scheme to work. This observation is
important since it allows the reduction scheme to be applied to the
class of mechanical systems with (possibly nonholonomic) kinematic
constraints; see \cite{B00,BvdS01} for further information and a
discussion of its relationship to other recent results in this area.
Of course, if $\D$ happens to be closed then also the reduced Dirac
structure $\D_\mu$ will be closed.

Finally, we would like to mention that in \cite{B00,BvdS01} it is
shown that the reduction scheme can also be applied the other way
around, starting with factorizing the symmetry group $G$ and
afterwards restricting the result to a level set of the remaining
first integrals (which actually turn out to be Casimir functions). The
resulting implicit Hamiltonian system on $\tilde{M}_\mu$ is
\emph{isomorphic} to the system $(M_\mu,\D_\mu,H_\mu)$. Notice that
$\tilde{M}_\mu$ is actually the \emph{orbit reduced space}
$P^{-1}(\Orb_\mu)/G$, where $\Orb_\mu$ denotes the coadjoint orbit in
$\mathfrak{g}^\ast$ through $\mu$. See also \cite{LM87,M92,O93} for
some classical references.

\paragraph{Intrinsic reductions} 
In this paragraph we will set the reduction results described above
against what we will call \emph{intrinsic reductions}. The latter are
independent of any symmetry properties of the Dirac structure and, in
fact, can be perfomed on \emph{any closed Dirac structure}. These kind
of reductions have been described in the literature by various authors
\cite{C90,D93,LWX98,B00,BvdS01}.

Consider a closed Dirac structure $\D$ on $M$. Then by condition
(\ref{eq:closed}) it follows that the characteristic distribution
$\Delta$ is involutive, and hence by Frobenius' Theorem defines a
regular foliation $\Phi_\Delta$ of $M$ into integral submanifolds of
$\Delta$. On the other hand, the distribution defined by\footnote{This
  distribution is usually denoted by $\mathsf{G}_1$ in the
  literature.}
\begin{equation}
\Theta(x) := \{ X(x) \mid  X \in \mathfrak{X}_{loc}(M), 
\exists \alpha \in \Omega_{loc}^1 (M), (X,\alpha) \in \mathfrak{D}_{loc} \}
\end{equation}
clearly is also involutive, defining a regular foliation $\Phi_\Theta$
of $M$ into integral submanifolds of $\Theta$. (For the moment we will
assume that both distributions are constant dimensional.)

There a two logical ways to ``reduce'' the Dirac structure on $M$ to a
lower dimensional manifold. The first is to project the Dirac
structure to the quotient manifold $M/\Phi_\Delta$, i.e., by factoring
out the characteristic distribution. This was done in \cite{C90} where
it was shown that the Dirac structure $\D$ on $M$ induces a well
defined Poisson bracket on the quotient manifold $M/\Phi_\Delta$
(\cite{C90}, Corollary~2.6.3). This remarkable result was generalized
in \cite{LWX98} to Dirac structures on Lie bialgebroids as described
in Remark \ref{rem:algebroids}, where is was refered to as
\emph{Poisson reduction}. In \cite{B00} it was observed that in fact
this reduction can be seen as a special case of \emph{symmetry
  reduction} if one notices that the distribution $\Delta$ is a
symmetry distribution of $\D$, i.e., every vector field $Y \in \Delta$
is a symmetry of $\D$ as in Definition \ref{def:sym}. The Dirac
structure $\D$ can be projected to a Dirac structure $\hat{\D}$ on
$M/\Phi_\Delta$ using Proposition \ref{prop:proj}. It turns out that
$\hat{\D}$ is exactly the graph of the Poisson structure corresponding
to the Poisson bracket defined by Courant \cite{C90}. We refer to
\cite{B00}, Example~4.2.4, p.~73, for more details.

The second reduction possibility is to restrict the Dirac structure to
each of the integral submanifolds of $\Theta$. This can be done using
Proposition \ref{prop:DN} and results in a Dirac structure on each of
the integral submanifolds of $\Theta$. In \cite{B00,BvdS01} it is
shown that each of the reduced Dirac structures represents a
presymplectic structure on the corresponding leaf of the foliation,
see \cite{B00}, Example~4.1.8, p.~69, and \cite{BvdS01}, Example~9, p.~79.
This corresponds to Theorem~2.3.6 in Courant \cite{C90} and
Theorem~2.2 in Dorfman \cite{D93}, stating that a closed Dirac
structure has a foliation by presymplectic leaves.

Once more we want to stress that the reductions described above are
``intrinsic'' and have nothing to do with the existence of any
symmetry groups of the implicit Hamiltonian system (although, as
explained above, the first reduction can be interpreted in terms of
symmetries of the Dirac structure). They can be perfomed on any
\emph{closed} Dirac structure. We will not concentrate on these
intrinsic reductions anymore, and instead will investigate
\emph{symmetry Lie groups} of implicit Hamiltonian systems, together
with their (singular) reductions. Doing so, we do \emph{not} assume
that the Dirac structure is closed and in fact all our results will be
presented for the general case.

\section{Singular reduction}\label{sec:sing}
\setcounter{equation}{0}

Contrary to the regular reduction reviewed in the previous section we
now describe in a purely \emph{topological} way how to obtain a
reduced Dirac structure on the reduced space $M_\mu$ if $M_\mu$ is not
a manifold. This is the case when $\mu$ is a singular value of the
momentum map $P$. In that case, vector fields and differential
one-forms on $M_\mu$ are not defined and therefore the results
described in the previous section cannot be used. Describing the
dynamics corresponding to such a topologically reduced Dirac structure
on $M_\mu$ will be done in section \ref{sec:sing_dyn}. For easiness of
exposition we will take $\mu = 0$ throughout the rest of this paper.

From this point on we specifically consider only Dirac structures of
the form (\ref{ds1}), admitting symmetries as described in Proposition
\ref{prop:sym}, that is, given a vector subbundle $\Delta \subset TM$,
and a generalized Poisson structure $J: T^\ast M \rightarrow TM$ (the
Jacobi identity does not necessarily hold), the Dirac structure is
defined by
\begin{equation}
\label{special Dirac structure}
\D(x) = \{ (v,v^\ast) \in T_xM \times T_x^\ast M \mid v - J(x) v^\ast \in \Delta(x), \; v^\ast \in \Delta^\circ (x) \}
\end{equation}
and $Y \in \mathfrak{X}(M)$ is a symmetry of $\D$ if 
\begin{equation}
\label{special symmetry}
L_Y J=0 \text{~and~} L_Y Z \text{~is a local section of~} \Delta \text{~ whenever~} Z \text{~is a local section of~} \Delta.
\end{equation}

Consider such a Dirac structure $\D$ on a manifold $M$ admitting a
symmetry Lie group $G$ with corresponding $Ad^\ast$-equivariant
momentum map $P$ satisfying (\ref{momentum}); it is not assumed that
$G$ acts regularly on $M$. Let $\mu =0 \in \mathfrak{g}^\ast$ be a
singular value of $P$ and consider the level set $N = P^{-1}(0)$ which
is not a smooth submanifold of $M$. However, $N$ is a closed subset of
$M$ and is a topological space relative to the induced subspace
topology. The level set $N$ is $G$-invariant so one can endow the
orbit space $M_0:= N/G = P^{-1}(0)/G$ with the quotient topology.
Denote by $\pi :N \rightarrow M_0$ the canonical projection map, that
is, $\pi$ maps $x \in N$ onto its orbit $G\cdot x \in M_0$.

Define the set of smooth functions on $M_0$ as follows (see \cite{CB97}, or for the original source see \cite{Schwarz}).

\begin{definition} 
\label{smooth functions on reduced space}
A continuous function $f_0$ on $M_0$ is called \emph{smooth}, denoted
by $f_0 \in C^\infty(M_0)$, if there exists a smooth $G$-invariant
function $f \in C^\infty(M)^G$ such that $f_0 \circ \pi = f
\vert_{P^{-1}(0)}$.
\end{definition} 

Given the singular reduced space $M_0$ together with its
topology and a set of smooth functions $C^\infty(M_0)$ on $M_0$, we
want to define a reduced Dirac structure on $M_0$. Let $\D$ be a Dirac
structure  on $M$ of the type defined in (\ref{special Dirac
  structure}) and assume that the infinitesimal generators $\xi_M, \;
\xi \in \mathfrak{g}$, satisfy the conditions in (\ref{special
  symmetry}). Since $G$ is a symmetry Lie group of the generalized
Poisson bracket $\{ \cdot,\cdot \}: C^\infty(M) \times C^\infty(M)
\rightarrow C^\infty(M)$, corresponding to the bundle map $J$, we can
use the theory in \cite{ACG91,CB97} to define a generalized Poisson
bracket $\{ \cdot,\cdot \}_0: C^\infty(M_0) \times C^\infty(M_0)
\rightarrow C^\infty(M_0)$ on the singular reduced space $M_0$. This
goes as follows. Let $f_0,h_0 \in C^\infty(M_0)$ and let $f,h \in
C^\infty(M)^G$ be such that $f_0 \circ \pi = f \vert_{P^{-1}(0)}$ and
$h_0 \circ \pi = h \vert_{P^{-1}(0)}$. Define the singular
reduced generalized bracket by
\begin{equation}\label{redbrack}
\{ f_0,h_0 \}_0 \circ \pi = \{ f,h \} \vert_{P^{-1}(0)}.
\end{equation}
This gives a well defined generalized Poisson bracket on $M_0$. In
particular, \eqref{redbrack} does not depend on the choice of the
$G$-invariant extensions $f$ and $h$ (whose existence is assumed, by
definition).

\begin{remark}{\rm
    The reduction theory in \cite{ACG91,CB97} is only developed for
    the singular reduction of \emph{symplectic} manifolds under a
    symmetry Lie group action. That is, the Poisson bracket $\{
    \cdot,\cdot \}$ is assumed to be nondegenerate and to satisfy the
    Jacobi identity. In principle, however, these results generalize
    immediately to the case of singular reduction of generalized
    Poisson brackets, as described above by (\ref{redbrack}). In
    particular, \cite{ACG91,CB97} show that (under the assumption that
    $G$ acts properly) nondegeneracy of $\{ \cdot,\cdot \}$ implies
    that of $\{ \cdot,\cdot \}_0$. Also, from (\ref{redbrack}) it
    follows immediately that $\{ \cdot,\cdot \}_0$ satisfies the
    Jacobi identity if $\{ \cdot,\cdot \}$ does.
    
    Once again for clarity: In this paper we do neither assume that
    the generalized Poisson bracket $\{ \cdot,\cdot \}$ is
    nondegenerate, nor that it satisfies the Jacobi identity.
    Furthermore, properness of the group action is not assumed until
    Section \ref{sec:orbitdecomp}.  
}\quad $\blacklozenge$
\end{remark}

Next, consider the vector subbundle $\Delta \subset TM$, defining
a constant dimensional distribution on $M$. We show that $\Delta$
defines a vector space $\hat{\Delta}$ consisting of derivations on the
space $C^\infty(M_0)$ of smooth function on $M_0$. Denote by
$\Gamma_{loc}(\Delta)$ the local sections of the subbundle $\Delta
\subset TM$. We show that every vector field $X \in
\Gamma_{loc}(\Delta)$ is ``tangent" to $N$. In the regular case, when
$N = P^{-1}(\mu)$ is a smooth submanifold of $M$, this means that $X$
restricts to a well defined vector field $\bar{X}$ on $N$. However, if
$\mu=0$ is a singular value of the momentum map, then $N$ is not a
smooth manifold and hence we have to define what ``tangent" means.
Recall that a vector field $X \in \mathfrak{X}(M)$ is in one-to-one
correspondence with a derivation, also denoted by $X: C^\infty(M)
\rightarrow C^\infty(M)$, on the set of smooth functions on $M$. The
correspondence is given by the formula\footnote{The derivation is
  usually called \emph{Lie derivative} and is also denoted by $X[f] =
  L_X f$.}
\begin{equation}
X[f] = \langle df, X \rangle, \; \forall f \in C^\infty(M).
\end{equation}
A derivation $X$ on $M$ is said to be \textit{tangent\/} to $N$ if it
restricts to a well defined derivation $\bar{X}$ on the set of Whitney
smooth functions on $N$. A continuous function $\bar{f}$ on $N$ is
said to be a \textit{Whitney smooth function\/} if there exists a
smooth function $f$ on $M$ such that $\bar{f} = f\vert_N$; the set of
Whitney smooth functions on $N$ is denoted by $W^\infty(N)$. Otherwise
stated, $X$ \textit{is tangent to $N$ if there exists a derivation
  $\bar{X}$ on $W^\infty(N)$ such that $X[f](x) = \bar{X}[f \vert_N]
  (x), \; \forall x \in N,$ for all $f \in C^\infty(M)$\/}. A
necessary and sufficient condition for $X$ to be tangent to $N$ is
that
\begin{equation}
X[f](x) = X[h](x), \; \forall x \in N,
\end{equation}
for all $f,h \in C^\infty(M)$ such that $f \vert_N = h\vert_N$. Notice
that in case $N$ is a smooth submanifold of $M$ and $N$ is
\emph{closed} in $M$, then the set $W^\infty(N)$ of Whitney smooth
functions on $N$ is equal to the set $C^\infty(N)$ of all smooth
functions on $N$ (as defined by the differential structure on the
submanifold $N$), and the above given definition yields the usual
meaning of a vector field being tangent to the submanifold $N$ (and
consequently restricting to a vector field on $N$).  \medskip

Consider a vector field (or equivalently, derivation) $X$ on $M$, and
define $\gamma(t)$ to be an \emph{integral curve} of $X$ through $x_0
\in M$ if\footnote{Take the coordinate functions $f=x^i$ to obtain the
  usual definition $\dot{\gamma}^i(t) = X^i(\gamma(t))$.}
\begin{equation}\label{curve}
\frac{d}{dt} f(\gamma(t)) = X[f](\gamma(t)), \quad \forall t, \; \forall f \in C^\infty(M), \; \gamma(0)=x_0.
\end{equation}
Now, let $X \in \Gamma_{loc}(\Delta)$ and $\gamma(t)$ be an integral curve of $X$ through $x_0 \in P^{-1}(0)$. In particular,
\begin{equation}
\frac{d}{dt} P_\xi(\gamma(t)) = X[P_\xi] (\gamma(t)) = 0, \; \forall t, \; \forall \xi \in \G,
\end{equation}
since by (\ref{ds1}) (or (\ref{special Dirac structure})) and
(\ref{momentum}), $dP_\xi(x) \in \Delta^\circ(x), \; \forall x \in M$.
This implies that the integral curve of $X \in \Gamma_{loc}(\Delta)$
through every $x_0 \in P^{-1}(0)$ is contained in $P^{-1}(0)$
(conservation of the momentum map). Then by the equivalence of
derivations and velocity vectors (remember that $M$ is a smooth
manifold) it follows that
\begin{equation}
X[f](x_0) = \frac{d}{dt} f(\gamma(t)) \vert_{t=0} = \frac{d}{dt} h(\gamma(t)) \vert_{t=0} = X[h](x_0),
\end{equation}
for all $f,h \in C^\infty(M)$ such that $f \vert_N = h \vert_N$. So we
have shown that every vector field $X \in \Gamma_{loc}(\Delta)$ is
tangent to $N= P^{-1}(0)$ and consequently restricts to a well defined
derivation $\bar{X}$ on $W^\infty(N)$. In conclusion, the constant
dimensional distribution $\Delta$ on $M$ restricts to a vector space
$\bar{\Delta}$ of derivations on $W^\infty(N)$. If $\Delta$ is locally
spanned by the independent vector fields $X_1,\dots,X_m$, then
$\bar{\Delta}$ is locally spanned by the independent derivations
$\bar{X}_1,\dots,\bar{X}_m$.

Using the above results we are able to show that the distribution
$\Delta$ on $M$ projects to a well defined vector space $\hat{\Delta}$
of derivations on the smooth functions $C^\infty(M_0)$. A vector field
$X$ on $M$ is said to \textit{project\/} to $M_0$ if there exists a
derivation $\hat{X}$ on $C^\infty(M_0)$ such that for every $f \in
C^\infty(M)^G$, $X[f](x) = \hat{X}[f_0](\pi(x)), \;\forall x \in N$,
with $f_0$ defined by $f_0 \circ \pi = f \vert_N$. It is clear that
$X$ restricts to a well defined derivation $\hat{X}$ on
$C^\infty(M_0)$ if and only if $X[f](x)$ does not depend on the
extension of $f_0 \circ \pi$ off $N$ to $M$ and furthermore $X[f](x) =
X[f](y)$ for all $x,y \in N$ such that $\pi(x)=\pi(y)$. Now let $X$ be
a local section of $\Delta$. Since $X$ is tangent to $N$ it follows
that $X[f](x) = \bar{X}[f \vert_N](x) = \bar{X}[f_0 \circ \pi](x), \;
\forall x \in N,$ and therefore its value does not depend on the
extension of $f_0 \circ \pi$ off $N$ to $M$. It remains to show that
\begin{equation}\label{projG_0}
X[f](x) = X[f](y),\; \forall x,y \in N \text{~such~that~} \pi(x) =\pi(y).
\end{equation}
In general, this is not true for every local section $X$ of $\Delta$.
However, we will show that \textit{there exists a basis of local
  sections $X_1, \dots,X_m$ of $\Delta$ which satisfies\/}
(\ref{projG_0}).

Denote by $\mathfrak{V}_{loc}$ the space of local sections of the
vertical distribution $V$ defined by $V(x):= \mbox{span}\:\{\xi_M(x)
\mid \xi \in \mathfrak{g} \}$. Since $L_{\xi_M} \Gamma_{loc}(\Delta)
\subset \Gamma_{loc}(\Delta)$ for every $\xi \in \mathfrak{g}$ by
(\ref{special symmetry}), it follows that
$[\Gamma_{loc}(\Delta),\mathfrak{V}_{loc}] \subset \mathfrak{V}_{loc}
+ \Gamma_{loc}(\Delta)$. Indeed, taking an arbitrary local section of
$V$ of the form $Y = \sum_i h_i \xi_M^i, \; h_i\in C^\infty(M),\;
i=1,\dots,r=\mbox{dim}\: \mathfrak{g}$, $\xi_M^1,\dots,\xi_M^r$ being
a local basis of $V$, and letting $X \in \Gamma_{loc}(\Delta)$, it
follows that
\begin{equation*}
[X,Y] = [X, \sum_{i=1}^r h_i \xi_M^i] = \sum_{i=1}^r h_i [X,\xi_M^i] + (L_Xh_i)  \xi_M^i \in \Gamma_{loc}(\Delta) + \mathfrak{V}_{loc}
\end{equation*}
which proves the inclusion $[\Gamma_{loc}(\Delta),\mathfrak{V}_{loc}]
\subset \mathfrak{V}_{loc} + \Gamma_{loc}(\Delta).$

Assuming that the distribution $V + \Delta$ has constant dimension on
$M$, the above inclusion implies that there exists a basis $X_1,
\dots,X_m$ of local sections of $\Delta$ such that
$[X_i,\mathfrak{V}_{loc}] \subset \mathfrak{V}_{loc}, \; i=1,\dots,m$;
see e.g.~Theorem 7.5 on page 214 in the book by Nijmeijer and van der
Schaft \cite{ns} (in the notation of that theorem: the involutive
distribution $D$ is $V$, the distribution $G$ is $\Delta$, and one
takes $f = 0$). In particular, $[X_i, \xi_M] \in \mathfrak{V}_{loc}$,
which implies that for all $f \in C^\infty(M)^G$
\begin{equation}
0 = [X_i, \xi_M][f] = X_i [L_{\xi_M}f] - L_{\xi_M}( X_i[f]) = - L_{\xi_M}( X_i[f]), \; \forall \xi \in \mathfrak{g}.
\end{equation}
This means that the function $X_i[f]$ is $G$-invariant and therefore
satisfies (\ref{projG_0}). In conclusion, there exists a basis
$X_1,\dots,X_m$ of local sections of $\Delta$ such that each $X_i$
projects to a well defined derivation $\hat{X}_i$ on $C^\infty(M_0)$.
The derivations $\hat{X}_1,\dots,\hat{X}_m$ locally span (in other
words, form a basis of) a vector space of derivations on
$C^\infty(M_0)$, denoted by $\hat{\Delta}$.

\begin{remark}\label{rem:proj}
  {\rm In the regular case, i.e., when $\mu=0$ is a regular value of
    the momentum map and $G$ acts freely and properly on $M$, the
    reduced space $M_0$ is a smooth manifold. Furthermore, the set of
    smooth functions $C^\infty(M_0)$ equals the set of smooth
    functions as defined by the differential structure on $M_0$.
    Indeed, since $N=P^{-1}(0)$ is closed in $M$, every smooth
    $G$-invariant function on $N$ can be smoothly extended to a
    $G$-invariant function on $M$ \cite{ACG91}. In that case the
    notion of a ``projecting derivation" as defined above has the
    usual meaning of projection of a vector field on $M$ to a vector
    field on the reduced space $M_0$. In particular there exists a
    basis $X_1,\dots,X_m$ of local sections of $\Delta$, which are
    tangent to $N$, such that the restrictions
    $\bar{X}_1,\dots,\bar{X}_m$ project to $M_0$. That is, each
    $\bar{X}_i$ is $\pi$-related to a vector field $\hat{X}_i$ on
    $M_0$. The projected vector fields $\hat{X}_1,\dots,\hat{X}_m$
    form a basis of local sections of $\hat{\Delta}$.  
} \quad $\blacklozenge$
\end{remark}

So far we have defined a set of smooth functions $C^\infty(M_0)$ on
the singular reduced space $M_0$, together with a generalized Poisson
bracket $\{ \cdot,\cdot \}_0$ and a vector space $\hat{\Delta}$ of
derivations on $C^\infty(M_0)$. Recall that the original Dirac
structure $\D$ (of the type given by (\ref{ds1}) or (\ref{special
  Dirac structure})) on the manifold $M$ was completely determined by
the generalized Poisson bracket $\{ \cdot,\cdot \}$, corresponding to
$J$, and the distribution $\Delta$. Therefore it makes sense to define
a reduced Dirac structure on $M_0$ as follows:

\begin{definition}\label{singds}
  Consider the singular reduced space $M_0$ together with the set of
  smooth functions $C^\infty(M_0)$. The \emph{singular reduced Dirac
    structure} $\D_0$ is defined as the pair $(\{ \cdot,\cdot \}_0,
  \hat{\Delta})$.
\end{definition}

We also call $\D_0$ a \emph{topological Dirac structure}. It will be
shown in the next section that the singular reduced Dirac structure
$\D_0$ defines a Hamiltonian formalism on the singular reduced space
$M_0$.

In order to better comprehend the meaning of the singular reduced
Dirac structure defined in Definition \ref{singds}, we show that in
the case of regular reduction the topological Dirac structure $\D_0$
exactly defines the regular reduced Dirac structure on $M_0$.

\paragraph{Regular reduction.}
Suppose that $\mu =0$ is a regular value of the momentum map and $G$
acts regularly on $M$ (that is, $M/G$ is a smooth manifold and the projection
$M \rightarrow M/G$ is a surjective submersion; for example, if $G$ acts
freely and properly, these conditions are satisfied). According to the
results described in Section \ref{sec:sr} the Dirac structure $\D$ on
$M$ is reduced to a Dirac structure $\hat{\D}$ on the manifold $M_0$
in two steps: firstly, $\D$ is restricted to a Dirac structure $\D_N$ on
$N=P^{-1}(0)$ defined by (\ref{DNsec}) and, secondly, $\D_N$ is
projected to a Dirac structure $\hat{\D}$ on $M_0$ defined by
(\ref{Dproj}) (with $\mathfrak{D}_{loc}$ replaced by
$(\mathfrak{D}_N)_{loc}$). Otherwise stated, in terms of its local
sections,
\begin{align}\nonumber
\hat{\mathfrak{D}}_{loc} = \{ (\hat{X},\hat{\alpha}) \in \mathfrak{X}_{loc}(M_0) \oplus \Omega^1_{loc} (M_0) \mid \; & \exists (X,\alpha) \in \mathfrak{D}_{loc} \;\mbox{such that $X$ is tangent to $N$}\\ 
&\mbox{and}\; \iota_\ast X \sim_\pi \hat{X}, \; \iota^\ast \alpha = \pi^\ast \hat{\alpha} \},
\end{align}
where $\iota: N \hookrightarrow M$ is the inclusion and $\iota_\ast X$
denotes the push forward of $X$ to a vector field on $N$ (that is,
$\iota_\ast X$ is simply the restriction of $X$ to $N$ which is
possible since $N$ is a closed submanifold of $M$ and $X$ is tangent
to $N$ by hypothesis).

Consider the topological Dirac structure $\D_0$ given by Definition
\ref{singds}. Since $M_0$ is a manifold, the generalized Poisson
bracket $\{ \cdot,\cdot \}_0$ on the set of smooth functions
$C^\infty(M_0)$ (see also Remark \ref{rem:proj}) defines a
skew-symmetric vector bundle map $J_0: T^\ast M_0 \rightarrow TM_0$ by
$J_0(df_0,dh_0) = \{ f_0,h_0 \}_0, \; f_0,h_0 \in C^\infty(M_0)$. The
vector space $\hat{\Delta}$ of derivations on $C^\infty(M_0)$ defines
a constant dimensional distribution of vector fields on $M_0$ (in
other words, a vector subbundle of $TM_0$), also denoted by
$\hat{\Delta}$. Then the topological Dirac structure $\D_0$ defines a
Dirac structure on the manifold $M_0$, also denoted by $\D_0$, which
in terms of its local sections is given by
\begin{equation}\label{D0}
(\mathfrak{D}_0)_{loc} = \{ (\hat{X},\hat{\alpha}) \in \mathfrak{X}_{loc}(M_0) \oplus \Omega^1_{loc}(M_0) \mid \hat{X} - J_0 \hat{\alpha} \in \Gamma_{loc}(\hat{\Delta}), \; \hat{\alpha} \in \Gamma_{loc}(\hat{\Delta}^\circ) \}.
\end{equation}
Indeed, $\D_0$ is a Dirac structure on $M_0$ as defined in Definition
\ref{def:ds} (notice that it is of the same form as in (\ref{ds2})).
We show that $\hat{\D} = \D_0$. Since both are Dirac structures and
therefore their fibers are of the same dimension (i.e., dim $M_0$), it
is enough to show that $\hat{\D} \subset \D_0$.

If $(\hat{X},\hat{\alpha})$ is a local section of $\hat{\D}$, then
there exists a local section $(X,\alpha)$ of $\D$ such that
$\iota_\ast X \sim_\pi \hat{X}$ and $\iota^\ast \alpha = \pi^\ast
\hat{\alpha}$. Since $(X,\alpha)$ is a local section of $\D$ one has
\begin{equation}\label{inD}
Z:= X - J\alpha \;\text{~is~a~local~section~ of~} \Delta, \; \; \alpha \text{~is~a~local~section~ of~} \Delta^\circ.
\end{equation}
Consider the vector field $J\alpha \in \mathfrak{X}_{loc}(M)$. Since
$(J\alpha,\alpha) \in \mathfrak{D}_{loc}$ it follows from
(\ref{momentum}) and $\D = \D^\perp$ that
\begin{equation}
(J\alpha)[P_\xi](x) = \langle dP_\xi, J\alpha \rangle(x) = - \langle \alpha, \xi_M \rangle(x) = - \langle \hat{\alpha}, 0 \rangle(\pi(x)) = 0, \; \forall x \in N,\;\forall \xi \in \mathfrak{g}.
\end{equation}
This implies that the vector field $J\alpha$ is tangent to $N$.
Furthermore, by construction of the reduced generalized bracket
(\ref{redbrack}) it follows that $\iota_\ast (J\alpha) \sim_\pi J_0
\hat{\alpha}$. Since also $\iota_\ast X \sim_\pi \hat{X}$, equation
(\ref{inD}) implies that there exists a vector field $\hat{Z} \in
\mathfrak{X}_{loc}(M_0)$ such that $\iota_\ast Z \sim_\pi \hat{Z}$. It
follows that $\hat{Z} \in \Gamma_{loc}(\hat{\Delta})$ by construction
of $\hat{\Delta}$. This yields
\begin{equation}
\hat{X} - J_0\hat{\alpha} = \hat{Z} \in \Gamma_{loc}(\hat{\Delta}).
\end{equation}
By construction, the distribution $\hat{\Delta}$ is spanned by vector
fields $\hat{Z}_1,\dots,\hat{Z}_m$ for which there exists a basis of
vector fields $Z_1,\dots,Z_m \in \Gamma_{loc}(\Delta)$ such that
$\iota_\ast Z_j \sim_\pi \hat{Z}_j, \;j=1,\dots,m$. Since $\iota^\ast
\alpha = \pi^\ast \hat{\alpha}$ and $\alpha \in
\Gamma_{loc}(\Delta^\circ)$, it follows immediately that
\begin{equation}
\langle \hat{\alpha},\hat{Z}_j \rangle \circ \pi = \langle \alpha,Z_j \rangle \circ \iota = 0, \; j=1,\dots,m,
\end{equation}
and therefore $\hat{\alpha} \in \Gamma_{loc}(\hat{\Delta}^\circ)$. In
conclusion, $(\hat{X},\hat{\alpha})$ is a local section of $\D_0$. So
we have shown that $\hat{\D} \subset \D_0$ and since both are Dirac
structures on $M_0$ this implies that $\hat{\D} = \D_0$.

We conclude that in the case of regular reduction the topological
Dirac structure $\D_0$ exactly defines the regular reduced Dirac
structure on $M_0$.

\section{Singular dynamics}\label{sec:sing_dyn}
\setcounter{equation}{0}

In this section a Hamiltonian formalism is described corresponding to
the singular reduced Dirac structure $\D_0$ of Definition
\ref{singds}. This formalism defines the dynamics corresponding to an
implicit Hamiltonian system $(M_0,\D_0,H_0)$ on the topological space
$M_0$. We show that if $(M_0,\D_0,H_0)$ is the reduction of the
implicit Hamiltonian system $(M,\D,H)$ to $M_0$, then the
$G$-projectable solutions of $(M,\D,H)$ project to solutions of the
reduced system $(M_0,\D_0,H_0)$.

First let us define a Hamiltonian formalism on a topological space
in the spirit of Sikorski differential spaces (see \cite{S, CSn91}).
Consider a topological space $M_0$ together with a subalgebra
$C^\infty(M_0)$ of the continuous functions on $M_0$, called the set
of smooth functions on $M_0$. A continuous curve $\gamma(t)$ on $M_0$
is said to be \emph{smooth} (see \cite{S}) if $f_0 \circ \gamma$ is smooth, as a function
from (a subinterval of) $\R$ to $\R$, for every $f_0 \in
C^\infty(M_0)$. Let $\hat{X}$ denote a derivation on $C^\infty(M_0)$.
An \emph{integral curve} of $\hat{X}$ through some point $x_0 \in M_0$
is defined (see \cite{S}) as a smooth curve $\gamma(t)$ for which, cf.~(\ref{curve}),
\begin{equation}
\frac{d}{dt} f_0(\gamma(t)) = \hat{X}[f_0](\gamma(t)), \quad \forall t, \; \forall f_0 \in C^\infty(M_0), \; \gamma(0)=x_0.
\end{equation}
Let $\D_0$ be a topological Dirac structure on $M_0$, consisting of a
generalized Poisson bracket $\{ \cdot,\cdot \}_0: C^\infty(M_0) \times
C^\infty(M_0) \rightarrow C^\infty(M_0)$ (that is, the Jacobi identity
does not necessarily hold) and a vector space $\hat{\Delta}$ of
derivations on $C^\infty(M_0)$. Furthermore let $H_0 \in
C^\infty(M_0)$ be a smooth function on $M_0$, called the Hamiltonian
function. Notice that $\{ \cdot, H_0 \}_0: C^\infty(M_0) \rightarrow
C^\infty(M_0)$ defines a derivation on $C^\infty(M_0)$ by $\{
\cdot,H_0 \}_0[f_0] : = \{ f_0,H_0 \}_0, \; f_0 \in C^\infty(M_0)$.
Furthermore, if $\hat{X}$ is a derivation on $C^\infty(M_0)$ and $x
\in M_0$, then $\hat{X}(x): C^\infty(M_0) \rightarrow \R$ is defined
by $(\hat{X}(x))[f_0] := \hat{X}[f_0](x), \; f_0 \in C^\infty(M_0)$.
The three-tuple $(M_0,\D_0,H_0)$ defines an implicit Hamiltonian
system in the following way:

\begin{definition}\label{def:sihs}
  A smooth curve $\gamma(t)$ on $M_0$ is called an \emph{integral
    curve} (or, \emph{solution}) of \linebreak[4]$(M_0,\D_0,H_0)$ if
  there exists a derivation $\hat{X}$ on $C^\infty(M_0)$ such that
  $\gamma(t)$ is an integral curve of $\hat{X}$ and
\begin{align}\label{sihs:1}
\hat{X}(&\gamma(t)) - \{ \cdot, H_0 \}_0(\gamma(t)) \in \hat{\Delta}(\gamma(t)), \;\forall t,\\[2ex]
\label{sihs:2}
&\hat{Z}[H_0](\gamma(t)) = 0, \; \forall t, \forall \hat{Z} \in \hat{\Delta}.
\end{align}
The implicit Hamiltonian system $(M_0,\D_0,H_0)$ is defined as the
total set of integral curves $\gamma(t)$ of $(M_0,\D_0,H_0)$.
\end{definition}

If $M_0$ is a smooth manifold, Definition
\ref{def:sihs} of an implicit Hamiltonian system equals Definition
\ref{def:ihs} given in Section \ref{sec:ihs} (with $\D_0$ defined by
(\ref{D0})). However, since in general $M_0$ is not a smooth manifold
but only a topological space, the implicit Hamiltonian system
$(M_0,\D_0,H_0)$ \emph{cannot} be written as a set of differential and
algebraic equations. As for implicit Hamiltonian systems defined on
manifolds, the implicit Hamiltonian system $(M_0,\D_0,H_0)$ is energy
conserving, cf.~(\ref{ec}),
\begin{equation}
\frac{dH_0}{dt}(\gamma(t)) = \hat{X}(H_0)(\gamma(t)) = \{ H_0,H_0 \}_0(\gamma(t)) = 0, \; \forall t.
\end{equation}

\begin{remark}
{\rm 
Equation (\ref{sihs:1}) implies that 
\begin{equation}\label{sihs:1rem}
\frac{d}{dt} f_0 (\gamma(t)) = \{ f_0, H_0 \}_0(\gamma(t)), \; \forall t, \forall f_0 \in \A_{\D_0},
\end{equation}
where $\A_{\D_0} = \{ f_0 \in C^\infty(M_0) \mid \hat{Z}[f_0] = 0, \;
\forall \hat{Z} \in \hat{\Delta} \}$. However, (\ref{sihs:1rem}) does
not imply (\ref{sihs:1}). Even in the regular case it is not true that
$\hat{Z}$ being a local section of $\hat{\Delta}$ is equivalent to
$\hat{Z}[f_0] = 0, \; \forall f_0 \in \A_{\D_0}$. A counterexample can
easily be constructed (by considering a suitable \emph{noninvolutive}
distribution $\hat{\Delta}$).  } \quad $\blacklozenge$
\end{remark}

\begin{remark}\label{rem:explicit}
  {\rm If $\hat{\Delta} = 0$ (i.e., $\Delta=0$, that is, no
    constraints), then $\A_{\D_0} = C^\infty(M_0)$ and
    (\ref{sihs:1},\:\ref{sihs:2}) are equivalent to (\ref{sihs:1rem}).
    In this case, the Hamiltonian dynamics defined by
    (\ref{sihs:1rem}) is exactly the singular reduced Hamiltonian
    dynamics as defined in \cite{CS91,CB97,OR98,OR03,SL91}.  
} \quad $\blacklozenge$
\end{remark}

Recall that implicit Hamiltonian systems defined on manifolds define a
set of differential and algebraic equations and, as a consequence, the
standard results on existence and uniqueness of solutions for
ordinary differential equations do not apply. As explained in Section
\ref{sec:ihs}, in general one cannot expect neither global existence
nor uniqueness of solutions of these systems. Therefore one cannot
expect global existence and uniqueness of solutions of implicit
Hamiltonian systems on topological spaces as defined in Definition
\ref{def:sihs}. In particular, all solutions necessarily lie in the
\emph{constraint space}
\begin{equation}
M_0^c = \{ x \in M_0 \mid \hat{Z}[H_0](x) = 0, \; \forall \hat{Z} \in \hat{\Delta} \}
\end{equation}
(a topological space whose topology is induced from $M_0$).
However, if $(M_0,\D_0,H_0)$ is the singular reduction of an implicit
Hamiltonian system $(M,\D,H)$, then we shall show next that the
$G$-projectable solutions of $(M,\D,H)$ (if they exist) project to
solutions of the reduced system $(M_0,\D_0,H_0)$.

Let $x(t)$ be a solution of $(M,\D,H)$ with $x(0)\in N=P^{-1}(0)$.
Then by (\ref{integral}), the curve $x(t)$ is contained in $N$. Now
assume that $x(t)$ is a $G$-\emph{projectable} solution, that is,
there exists a projectable derivation (i.e., vector field) $X$ on
$C^\infty(M)^G$, which projects to a well defined derivation $\hat{X}$
on $C^\infty(M_0)$, such that $x(t)$ is an integral curve of $X$
(i.e., $\dot{x}(t) = X(x(t))$). By (\ref{ds1}) and (\ref{def:ihs-sol})
\begin{align}\label{sol1}
X(x(t)) &- \{ \cdot,H \}(x(t)) =: Z(x(t)) \in \Delta(x(t)), \; \forall t,\\[2ex]
&Y[H](x(t)) = 0, \; \forall t, \; \forall Y \in \Gamma_{loc}(\Delta).
\end{align}
Let $M_0$ be the singular reduced space and $\D_0$ be the singular
reduced Dirac structure on $M_0$. Since $H$ is assumed to be
$G$-invariant, its restriction to $N$ projects to a well defined
function $H_0 \in C^\infty(M_0)$ defined by $H_0 \circ \pi =
H\vert_N$; for the definition of $C^\infty(M_0)$ see Definition
\ref{smooth functions on reduced space} in Section \ref{sec:sing}.
Define the singular reduced implicit Hamiltonian system
$(M_0,\D_0,H_0)$ as in Definition \ref{def:sihs}. Project the curve
$x(t)$ to $M_0$ to obtain the smooth curve $\gamma(t) =\pi(x(t))$ on
$M_0$. Then $\gamma(t)$ is an integral curve of the derivation
$\hat{X}$. Indeed, take an arbitrary $f_0 \in C^\infty(M_0)$ and let
$f \in C^\infty(M)^G$ be such that $f_0 \circ \pi = f\vert_N$. Then
\begin{equation}
\frac{d}{dt} f_0 (\gamma(t)) = \frac{d}{dt} f(x(t)) = X[f](x(t)) = \hat{X}[f_0](\gamma(t)),\; \forall t,
\end{equation}
where we used the fact that $x(t)$ is an integral curve of $X$,
cf.~(\ref{curve}), and that $X$ projects to a derivation $\hat{X}$ on
$C^\infty(M_0)$. Furthermore, if $Y_1, \dots,Y_m$ is a basis of
projectable local sections of $\Delta$, projecting to a basis
$\hat{Y}_1,\dots,\hat{Y}_m$ of $\hat{\Delta}$, then it follows that
\begin{equation}
0 = Y_j[H](x(t)) = \hat{Y}_j[H_0](\gamma(t)), \; \forall t, j=1,\dots,m,
\end{equation}
which yields equation (\ref{sihs:2}). It remains to be proved that
(\ref{sihs:1}) is satisfied. Notice that by (\ref{redbrack}) the
derivation $\{ \cdot,H\}$ projects to a well defined derivation $\{
\cdot,H_0\}_0$ on $C^\infty(M_0)$. Since also $X$ projects to a
derivation $\hat{X}$ on $C^\infty(M_0)$ it follows that the derivation
$Z$ projects to a well defined derivation $\hat{Z}$ on
$C^\infty(M_0)$. Since $Y_1,\dots,Y_m$ is a projectable \emph{basis}
of local sections of $\Delta$ it follows from (\ref{sol1}) that at
each point $x_0$ on the curve $x(t)$, one has
\begin{equation}
Z(x_0) = c_1 Y_1 (x_0) + \dots + c_m Y_m (x_0), 
\end{equation}
for some constants $c_1,\dots,c_m \in \R$. We claim that 
\begin{equation}\label{eqZ}
\hat{Z}(\gamma_0) = c_1 \hat{Y}_1(\gamma_0) + \dots + c_m \hat{Y}(\gamma_0), \; \gamma_0=\pi(x_0).
\end{equation}
Indeed, take an arbitrary $f_0 \in C^\infty(M_0)$ and let $f \in C^\infty(M)^G$ be such that $f_0 \circ \pi = f\vert_N$. Then
\begin{align}\nonumber
(\hat{Z}(\gamma_0))[f_0] = \hat{Z}[f_0](\gamma_0) &= Z[f](x_0)\\
\nonumber
&= \left(c_1 Y_1[f] + \dots + c_m Y_m[f]\right)(x_0)\\
\nonumber
&= \left(c_1 \hat{Y}_1[f_0] + \dots + c_m \hat{Y}_m[f_0]\right)(\gamma_0)\\
&= \left(c_1 \hat{Y}_1(\gamma_0) + \dots + c_m \hat{Y}_m(\gamma_0) \right) [f_0],
\end{align}
which yields (\ref{eqZ}). Since $\hat{Y_1},\dots,\hat{Y}_m$ forms a
basis of $\hat{\Delta}$ it follows that $\hat{Z}(\gamma(t)) \in
\hat{\Delta}(\gamma(t)), \; \forall t$. Therefore (\ref{sihs:1}) is
satisfied, which implies that $\gamma(t)$ is an integral curve of the
reduced implicit Hamiltonian system $(M_0,\D_0,H_0)$. We have obtained
the following

\begin{proposition}\label{prop:Gproj}
  Every $G$-projectable solution $x(t)$ of $(M,\D,H)$ with $x(0) \in
  P^{-1}(0)$ projects to a solution $\gamma(t) = \pi(x(t))$ of the
  singular reduced implicit Hamiltonian system $(M_0,\D_0,H_0)$.
\end{proposition}

\begin{remark}
  {\rm Suppose that the implicit Hamiltonian system $(M,\D,H)$ is of
    index 1. As remarked in Section \ref{sec:ihs}, the system can be
    restricted to an explicit Hamiltonian system on the constraint
    manifold $M_c$ defined by a generalized Poisson bracket denoted by
    $\{ \cdot,\cdot \}_c: C^\infty(M_c) \times C^\infty(M_c)
    \rightarrow C^\infty(M_c)$. The $G$-action leaves the manifold
    $M_c$ invariant and it follows that $G$ is a symmetry Lie group of
    the explicit Hamiltonian system on $M_c$, i.e., $L_{\xi_{M_c}} \{
    f,g \}_c = \{L_{\xi_{M_c}} f,g \}_c + \{ f, L_{\xi_{M_c}}g \}_c$,
    $\forall f, g \in C^\infty(M_c)$ and $L_{\xi_{M_c}}H_c = 0$ (where
    $H_c = H\vert_{M_c}$), for all $\xi \in \G$. The corresponding
    equivariant momentum map is given by the restriction of $P$ to
    $M_c$. Furthermore, every solution of the (restricted) system is
    $G$-projectable \cite{vdS98,B00,BvdS01}. We can use the singular
    reduction theory developed in \cite{CS91,CB97,OR98,SL91}, or
    equivalently the theory developed in this paper by considering
    $\Delta=0$, to reduce the system to a Hamiltonian system on the
    singular reduced space $(M_c)_0$. The reduced generalized Poisson
    bracket $(\{ \cdot,\cdot \}_c)_0$ is defined analogously to
    (\ref{redbrack}). The reduced dynamics is given by equation
    (\ref{sihs:1rem}), with $\{ \cdot,\cdot \}_0$ replaced by $(\{
    \cdot,\cdot \}_c)_0$ and $\A_{\D_0} = C^\infty((M_c)_0)$. Global
    existence of solutions now follows from Proposition
    \ref{prop:Gproj}. Furthermore, if the $G$-action is \emph{proper}
    then also uniqueness of solutions of the singular reduced system
    can be proved \cite{CB97,SL91}.  } \quad $\blacklozenge$
\end{remark}

\section{Orbit type decomposition}\label{sec:orbitdecomp}
\setcounter{equation}{0}

Consider a \emph{symplectic} manifold $(M,\omega)$ admitting a
symmetry Lie group with a corresponding momentum map. Let $M_0$ denote
the singular reduced space and $\{ \cdot,\cdot \}_0$ the singular
reduced Poisson bracket. The singular reduced Hamiltonian dynamics is
defined by equation (\ref{sihs:1rem}), cf.~Remark \ref{rem:explicit}.
In \cite{BL97,CB97,CS91,CSn91,OR03,SL91} it is shown that the space $M_0$
may be decomposed into a family of symplectic manifolds, called
pieces. The decomposition is by orbit type and defines a
stratification of the singular reduced space $M_0$. Furthermore, the
Hamiltonian flow corresponding to (\ref{sihs:1rem}) leaves the pieces
invariant and restricts to a regular Hamiltonian flow on each of the
pieces. In this section we show that these results can be generalized
to singular reduced implicit Hamiltonian systems. We will treat only
the special class of Dirac structures for which the generalized
Poisson bracket is nondegenerate.

Consider an implicit Hamiltonian system $(M,\D,H)$ with a Dirac
structure $\D$ as defined in (\ref{ds1}), or (\ref{special Dirac
  structure}), where the generalized Poisson bracket is assumed to be
nondegenerate and defined by a nondegenerate two-form $\omega$ on $M$,
i.e., the generalized Poisson structure $J: T^\ast M \rightarrow TM$
is given by $\omega^{-1}$, that is, $\{ f,h\} = \omega (X_f, X_h),\;
f,h \in C^\infty(M)$, where $X_f$ is defined by $df =
\omega(X_f,\cdot)$ and analogously for $X_h$. Notice that we do
\emph{not} assume that $\omega$ is a closed two-form. This means that
$\{ \cdot,\cdot \}$ does not necessarily satisfy the Jacobi identity.

Let $G$ be a symmetry Lie group of $(M,\D,H)$ as in Proposition
\ref{prop:sym}. From now on we will assume that the action of $G$ is
\emph{proper}. Assume that there exists an $Ad^\ast$-equivariant
momentum map $P: M \rightarrow \mathfrak{g}^\ast$ such that
\begin{equation}\label{mom}
dP_\xi = \omega (\xi_M,\cdot) \quad \mbox{(equivalently}\; \xi_M = JdP_\xi\mbox{) and}\quad dP_\xi \in \Gamma(\Delta^\circ), \; \forall \xi \in \mathfrak{g}.
\end{equation}
Notice that this implies, but is not equivalent to, (\ref{momentum}).

\begin{remark}\label{rem:mechsym}
  {\rm Symmetry groups as described above commonly occur within the
    class of constrained mechanical systems as described in Example
    \ref{ex:ms}. Consider the distribution $\Lambda = \ker A^T(q)$ on
    $Q$ and let $G$ be a Lie group acting properly on the
    configuration space $Q$ leaving $\Lambda$ invariant, that is,
    $L_{\xi_Q} \Gamma_{loc}(\Lambda) \subset \Gamma_{loc}(\Lambda), \;
    \forall \xi \in \mathfrak{g}$. The action of $G$ on $Q$ lifts to
    an action of $G$ on $M=T^\ast Q$ as follows (recall that
    $\omega=dp \wedge dq$ denotes the canonical symplectic form on
    $M$): define the infinitesimal generator $\xi_M$ to be the
    Hamiltonian vector field corresponding to the function $P_\xi(q,p)
    = p^T\xi_Q(q)$, i.e., $dP_\xi = \omega (\xi_M,\cdot), \; \forall
    \xi \in \mathfrak{g}$. Now recall that $\Delta$ is defined as the
    image of the matrix (\ref{mechDelta}). Then by construction,
    $L_{\xi_M} \omega = 0$ and $L_{\xi_M} \Gamma_{loc}(\Delta) \subset
    \Gamma_{loc}(\Delta), \; \forall \xi \in \mathfrak{g}$, and
    therefore $G$ defines a symmetry group of the implicit Hamiltonian
    system on $M$. The $Ad^\ast$-equivariant momentum map is defined
    by $P(q,p)(\xi) = p^T\xi_Q(q),\; \forall \xi \in \mathfrak{g}$,
    and hence satisfies the first condition in (\ref{mom}). If we
    assume furthermore that the symmetry group is ``horizontal", i.e.,
    $\xi_Q \in \ker A^T(q), \;\forall \xi \in \mathfrak{g}$, then also
    the second condition in (\ref{mom}) is satisfied. 
} \quad $\blacklozenge$ 
\end{remark}

The manifold $M$ can be decomposed into submanifolds as follows
\cite{CB97,OR03,SL91}. Let $K$ be a compact subgroup of $G$ and define
$M_{(K)}$ to be the set of points in $M$ whose stabilizer group $G_x=
\{g \in G \mid \phi(g,x)=x \}$ is conjugate to $K$, i.e.,
\begin{equation}
M_{(K)} = \{ x \in M \mid \exists g \in G \;\mbox{such that}\; gG_xg^{-1} = K \}.
\end{equation}
Notice that since the $G$-action is assumed to be proper, every
stabilizer group $G_x, \; x \in M$, is a compact subgroup of $G$.
$M_{(K)}$ is a submanifold of $M$ called the manifold of \emph{orbit
  type} $(K)$. On the set of compact subgroups of $G$ define an
equivalence relation by saying that $\tilde{K} \sim K$ if and only if
$\tilde{K}$ is conjugate to $K$. The equivalence class of $K$ is
denoted by $(K)$. As $(K)$ runs over the set of equivalence classes,
the manifolds $M_{(K)}$ partition $M$. Since the $G$-action is proper
this partition is locally finite. This is called the \emph{orbit type
  decomposition} of $M$.

Next we show that the image of the tangent of the momentum map at the
point $x \in M$ is equal to the annihilator in $\mathfrak{g}^\ast$ of
the Lie algebra of the stabilizer group $G_x$, i.e.,
\begin{equation}\label{annstab}
\mbox{Im}\; T_xP = \mathfrak{g}_x^\circ, \; \forall x \in M,
\end{equation}
see also \cite{AM78, CB97, LM87, MR99, OR03}. Indeed 
\begin{equation}
\xi \in \mathfrak{g}_x \Leftrightarrow \xi_M(x) = 0 \Leftrightarrow dP_\xi(x) = 0 \Leftrightarrow (T_xP( v))\xi = 0, \; \forall v \in T_xM \Leftrightarrow \xi \in (\mbox{Im}\; T_xP)^\circ,
\end{equation}
where we used (\ref{mom}) and the fact that $\omega$ is nondegenerate.
This yields that $\mathfrak{g}_x = (\mbox{Im}\; T_xP)^\circ$. Taking
the annihilator of both sides yields (\ref{annstab}).

Equation (\ref{annstab}) implies that the tangent of the restriction
of the momentum map $P$ to the manifold $M_{(K)}$ has constant rank
equal to the codimension of $K$ in $G$. It follows that the
intersection $P^{-1}(0) \cap M_{(K)}$ is a smooth submanifold of $M$.
Furthermore, the manifold $P^{-1}(0) \cap M_{(K)}$ is invariant under
the action of $G$. It turns out that the quotient
$(M_0)_{(K)}:=(P^{-1}(0) \cap M_{(K)})/G = \pi (P^{-1}(0) \cap
M_{(K)})$ is a smooth manifold \cite{CB97,OR03,SL91}. Consequently,
the singular reduced space $M_0$ is decomposed into a disjoint set of
manifolds, called pieces,
\begin{equation}\label{decomp}
M_0 = \coprod_{(K)} \;(M_0)_{(K)},
\end{equation}
where $(K)$ runs over the set of conjugacy classes of compact
subgroups of $G$. Since the orbit type decomposition of $M$ is locally
finite, the decomposition of $M_0$ is also locally finite.

Next let us define a generalized Poisson bracket on each of the
manifolds $(M_0)_{(K)}$. For clarity of exposition, consider the
following commuting diagram: \unitlength=5mm
\begin{center}
\begin{picture}(9,10)
\put(1,5){\makebox(0,0){$P ^{-1} (0)$}} 
\put(9,5){\makebox(0,0){$P ^{-1} (0)\cap M _{ (K)}$}}
\put(1,0){\makebox(0,0){$M_0$}} 
\put(9,0){\makebox(0,0){$(M_0)_{(K)}$}} 
\put(5,10){\makebox(0,0){$M$}}
\put(1,4){\vector(0,-1){3}}
\put(9,4){\vector(0,-1){3}}
\put(1,6){\vector(1,1){3}}
\put(9,6){\vector(-1,1){3}} 
\put(6.2,5){\vector(-1,0){3.7}}
\put(7.4,0){\vector(-1,0){5.6}}
\put(1,2){\makebox(0,0){$$}}
\put(10,2.5){\makebox(0,0){$\pi_{(K)}$}}
\put(0.4,2.5){\makebox(0,0){$\pi$}}
\put(1.6,8){\makebox(0,0){$\iota$}}
\put(8.4,8){\makebox(0,0){$\iota_{(K)}$}}
\put(5,5.6){\makebox(0,0){$\widetilde{ \iota}_{(K)}$}}
\put(5,.7){\makebox(0,0){$\iota^0_{(K)}$}}
\end{picture}
\end{center}
Here $\iota_{(K)}$ denotes the inclusion map and $\pi_{(K)}$ the
restriction of $\pi$ to $P^{-1}(0) \cap M_{(K)}$. The inclusions
$\tilde{\iota}_{(K)}$ and $\iota_{(K)}^0$ are defined such that the
diagram commutes. Define the set of Whitney smooth functions
$W^\infty((M_0)_{(K)})$ on $(M_0)_{(K)}$ as follows: a continuous
function $\bar{f}_0$ on $(M_0)_{(K)}$ is said to be a \textit{Whitney
  smooth function\/} if there exists a smooth $G$-invariant function
$f \in C^\infty(M)^G$ such that $\bar{f}_0 \circ \pi_{(K)} = f
\vert_{P^{-1}(0) \cap M_{(K)}}$. In fact, $W^\infty((M_0)_{(K)})$ is
equal to the set of functions obtained by the restriction of the
functions in $C^\infty(M_0)$ to $(M_0)_{(K)}$ (this is why it is
called the set of Whitney smooth functions). Indeed, the $G$-invariant
function $f$ descends to a smooth function $f_0$ on $M_0$, whose
restriction to $(M_0)_{(K)}$ is precisely $\bar{f}_0$. This can be
seen as follows:
\begin{equation}\label{f_restriction}
f_0 \circ \iota_{(K)}^0 \circ \pi_{(K)} = f_0 \circ \pi \circ \tilde{\iota}_{(K)} = f \circ \iota \circ \tilde{\iota}_{(K)} = f \circ \iota_{(K)} = \bar{f}_0 \circ \pi_{(K)},
\end{equation}
and since $\pi_{(K)}$ is surjective the result follows. Analogously to
(\ref{redbrack}) define a generalized Poisson bracket on
$W^\infty((M_0)_{(K)})$ by
\begin{equation}\label{redbrack(K)}
\{ \bar{f}_0,\bar{h}_0 \}_{(K)} \circ \pi_{(K)} = \{ f,h \} \vert_{P^{-1}(0)\cap M_{(K)}}.
\end{equation}
We need to show that the bracket is well defined and does not depend
on the choice of $G$-invariant extensions $f$ and $h$. If we can prove
that for every $f \in C^\infty(M)^G$ the flow of the Hamiltonian
vector field $X_f = \{ \cdot, f \}$ preserves the submanifold
$P^{-1}(0) \cap M_{(K)}$, then ${\cal I} = \{ f \in C^\infty(M)^G \mid
f\vert_{P^{-1}(0)\cap M_{(K)}} = 0 \}$ is a Poisson ideal of
$C^\infty(M)^G$ and it follows that the bracket is well defined, see
e.g.~\cite{CB97}, Appendix B, (5.2) \& (5.3). Now consider $f \in
C^\infty(M)^G$, that is, $\phi_g ^\ast f = f$, for all $g \in G$,
where $\phi_g: M \rightarrow M$ denotes the diffeomorphism given by
the action of the element $g \in G$ on $M$. Thus, since the action is
symplectic, it follows that $\phi_g^\ast X_f = X_{\phi_g^\ast f} =
X_f$ for all $g \in G$, which implies that the flow $\psi_t^f$ of
$X_f$ commutes with $\phi_g$ for every $g \in G$. We need to show that
$\psi_t^f$ preserves $M_{(K)}$. Since $M_{(K)} = G \cdot M_K$, where
$M_K = \{ x \in M \mid G_x = K \}$ (which is a submanifold of $M$) and the flow of $X_f$ and the
$G$-action commute, it is enough to show that $\psi_t^f$ preserves
$M_K$. Now, for every $g \in K$ and $x \in M_K$
\begin{equation}
\phi_g (\psi_t^f (x)) = \psi_t^f (\phi_g (x)) = \psi_t^f (x),
\end{equation}
since $\phi_g(x) = x$. Therefore, $K \subset G_{\psi_t^f(x)}$ (the
stabilizer group of $\psi_t^f(x)$). Suppose $g \in G_{\psi_t^f(x)}$,
i.e., $\phi_g (\psi_t^f (x))= \psi_t^f (x)$. Since the flow and the
$G$-action commute, this implies that $\psi_t^f (\phi_g (x)) =
\psi_t^f (x)$. Applying $\psi_{-t}^f$ to this relation yields $\phi_g
(x)= x$ so $g \in K$, i.e., $G_{\psi_t^f(x)} \subset K$. It follows
that $G_{\psi_t^f(x)} = K$ and therefore $\psi_t^f (x) \in M_K$. So we
have shown that for any $G$-invariant function $f$ the flow $\psi_t^f$
preserves the submanifold $P^{-1}(0) \cap M_{(K)}$. This implies that
the bracket (\ref{redbrack(K)}) is well defined.

\begin{remark}\label{rem:functions}
  {\rm The set of Whitney smooth functions $W^\infty((M_0)_{(K)})$ is
    dense in the set $C^\infty((M_0)_{(K)})$ of smooth functions as
    defined by the differential structure on $(M_0)_{(K)}$. Indeed,
    the pullback to $P^{-1}(0)\cap M_{(K)}$ of a smooth function
    $\bar{f}_0 \in C^\infty((M_0)_{(K)})$ compactly supported on
    $(M_0)_{(K)}$ can be extended to a smooth $G$-invariant function
    $f$ on $M$. Thus the bracket in (\ref{redbrack(K)})
    is a well defined generalized Poisson bracket $\{ \cdot,\cdot
    \}_{(K)}$ on $C^\infty((M_0)_{(K)})$. The associated generalized
    Poisson structure $J_{(K)}$ is the inverse of the nondegenerate
    two-form $\omega_{(K)}$ defined by the condition that the pullback
    of $\omega_{(K)}$ to $P^{-1}(0)\cap M_{(K)}$ equals the
    restriction of $\omega$ to $P^{-1}(0)\cap M_{(K)}$, see also
    \cite{CB97,OR03,SL91}.  } \quad $\blacklozenge$
\end{remark}

Having defined a generalized Poisson bracket on $(M_0)_{(K)}$ we now
show that the distribution $\Delta$ projects to a distribution
$(\hat{\Delta})_{(K)}$ on the manifold $(M_0)_{(K)}$. Recall from
Section \ref{sec:sing} that there exists a basis $X_1,\dots,X_m$ of
local sections of $\Delta$ such that $[\xi_M, X_j ] \in
\mathfrak{V}_{loc}, \; \forall \xi \in \mathfrak{g},\; j=1,\dots,m$.
By Proposition \ref{prop:sym} it follows that
\begin{equation}
[\xi_M, X_j] \in \Gamma_{loc}(\Delta \cap V), \quad \forall \xi \in \mathfrak{g}, \; j=1,\dots,m.
\end{equation}
If we make the assumption that $\Delta \cap V = 0$, then it follows
that $[\xi_M, X_j] =0, \; \forall \xi \in \mathfrak{g}, \;
j=1,\dots,m$. This means that the flow of $X_j$ commutes with the
$G$-action and therefore preserves the submanifold $P^{-1}(0) \cap
M_{(K)}$. This implies that $X_j$ is tangent to $P^{-1}(0) \cap
M_{(K)}$ and restricts to a local vector field $\bar{X}_j$ on
$P^{-1}(0) \cap M_{(K)}$, for all $j=1,\dots,m$. The basis
$\bar{X}_1,\dots,\bar{X}_m$ projects to a set of independent local
vector fields $\hat{X}_1,\dots,\hat{X}_m$ on $(M_0)_{(K)}$ which
define a distribution $(\hat{\Delta})_{(K)}$ on $(M_0)_{(K)}$.

\begin{remark}
{\rm 
If the system is of index 1, then it is always true that $\Delta(x) \cap V(x) = 0, \; \forall x \in M_c$ (see \cite{B00,BvdS01}).
} \quad $\blacklozenge$
\end{remark}

The generalized Poisson bracket $\{ \cdot,\cdot \}_{(K)}$, with its
generalized Poisson structure $J_{(K)}$, and the distribution
$(\hat{\Delta})_{(K)}$ define a Dirac structure $\D_{(K)}$ on
$(M_0)_{(K)}$, given in terms of its local sections by
\begin{align}\label{ds(K)}
(\mathfrak{D}_{(K)})_{loc} &= 
\left\{ (\hat{X},\hat{\alpha}) \in \mathfrak{X}_{loc}\left((M_0)_{(K)}\right) \times \Omega^1_{loc}\left((M_0)_{(K)}\right) \mid 
\right. \nonumber \\
&\qquad \left. \hat{X} - J_{(K)}\hat{\alpha} \in \Gamma_{loc}\left((\hat{\Delta})_{(K)}\right),\; \hat{\alpha} \in \Gamma_{loc}\left(\left((\hat{\Delta})_{(K)}\right)^\circ \right) \right\}.
\end{align}
The results in Section \ref{sec:sing} imply that $\D_{(K)}$ is exactly
the regular reduced Dirac structure on $(M_0)_{(K)}$. We summarize the
previous discussion in the following statement.

\begin{proposition}\label{prop:pieces}
  Consider an implicit Hamiltonian system $(M,\D,H)$ with a Dirac
  structure $\D$ of the type given in (\ref{ds1}), defined by a
  nondegenerate generalized Poisson bracket $\{ \cdot,\cdot \}$ and a
  smooth vector subbundle $\Delta$ of $TM$. Assume the system admits a
  symmetry Lie group $G$ acting properly on $M$ and satisfying the
  hypotheses of Proposition \ref{prop:sym}. Suppose also that the
  action admits an $Ad^\ast$-equivariant momentum map $P$ satisfying
  (\ref{mom}). Assume furthermore that $\Delta \cap V= 0$ and that the
  assumptions in Propositions \ref{prop:DN} and \ref{prop:proj} are
  satisfied.\footnote{Notice that $\Delta$ and $V$ being constant
    dimensional distributions and $\Delta \cap V =0$ imply that
    $\Delta + V$ is a constant dimensional distribution, i.e., a
    smooth vector subbundle of $TM$.} Then the singular reduced space
  $M_0$ is decomposed into a disjoint set of manifolds $(M_0)_{(K)}$,
  cf.~(\ref{decomp}), called pieces. The system $(M,\D,H)$ reduces, by
  means of regular reduction, to an implicit Hamiltonian system
  $((M_0)_{(K)},\D_{(K)},H_{(K)})$ on each of the pieces. Here
  $\D_{(K)}$ is the Dirac structure defined in (\ref{ds(K)}) and
  $H_{(K)} \in C^\infty((M_0)_{(K)})$ is defined by $H_{(K)} \circ
  \pi_{(K)} = H \vert_{P^{-1}(0) \cap M_{(K)}}$.
\end{proposition}

\begin{remark}
  {\rm The only reason why we assumed that the Dirac structure $\D$ is
    defined by a nondegenerate generalized Poisson bracket and the
    momentum map $P$ satisfies (\ref{mom}), is to show the equality
    (\ref{annstab}). This equality implies that the subsets $P^{-1}(0)
    \cap M_{(K)}$ are smooth manifolds and consequently the pieces
    $(M_0)_{(K)}$ are smooth manifolds as well. The construction of
    the reduced Dirac structure $\D_{(K)}$ however is completely
    general. This means that the result of Proposition
    \ref{prop:pieces} is valid for general Dirac structures and
    $Ad^\ast$-equivariant momentum maps of the form (\ref{momentum}),
    as long as it is certain that the pieces $(M_0)_{(K)}$ are smooth
    manifolds. The same will hold for the results in the remaining
    part of this section.\footnote{One could even think of further
      generalizing the results by allowing the pieces to be
      topological spaces instead of manifolds.}  
} \quad $\blacklozenge$
\end{remark}

Next we show that the regular reduced implicit Hamiltonian systems
$((M_0)_{(K)},\D_{(K)},H_{(K)})$ are precisely the restriction of the
singular reduced implicit Hamiltonian system $(M_0,\D_0,H_0)$ to the
pieces $(M_0)_{(K)}$. Indeed, first we show that the inclusion
$\iota_{(K)}^0 : (M_0)_{(K)} \rightarrow M_0$ is a Poisson map.
Consider $f_0,h_0 \in C^\infty(M_0)$ and their restrictions
$\bar{f}_0, \bar{h}_0 \in W^\infty((M_0)_{(K)})$ and let $f,h \in
C^\infty(M)^G$ be such that $f_0 \circ \pi = f\vert_{P^{-1}(0)}$, and
analogously for $h$. The definition of the brackets and the
commutativity of the diagram implies that (notice that by
(\ref{f_restriction}) $\bar{f}_0 \circ \pi_{(K)} = f\vert_{P^{-1}(0)
  \cap M_{(K)}}$ and analogously for $\bar{h}_0$)
\begin{equation}
\{ \bar{f}_0,\bar{h}_0 \}_{(K)} \circ \pi_{(K)} = \{ f,h \} \circ \iota_{(K)} = \{ f,h\} \circ \iota \circ \tilde{\iota}_{(K)} = \{ f_0,h_0 \}_0 \circ \pi \circ \tilde{\iota}_{(K)} = \{ f_0,h_0\}_0 \circ \iota_{(K)}^0 \circ \pi_{(K)}.
\end{equation}
Since $\pi_{(K)}$ is surjective it follows that 
\begin{equation}
\{ \bar{f}_0,\bar{h}_0 \}_{(K)} = \{ f_0,h_0\}_0 \vert_{(M_0)_{(K)}}.
\end{equation}
By construction, it is immediately clear that
$(\hat{\Delta})_{(K)}$ is the restriction of $\hat{\Delta}$ to
$W^\infty((M_0)_{(K)})$, defining by denseness a distribution on
$C^\infty((M_0)_{(K)})$. We conclude that the regular reduced Dirac
structure $\D_{(K)}$ is precisely the restriction of the singular
reduced Dirac structure $\D_0$ to the piece $(M_0)_{(K)}$. Next,
consider a solution $\gamma(t)$ of $(M_0,\D_0,H_0)$, as defined in
Definition \ref{def:sihs}. Since $\hat{X}_1,\ldots,\hat{X}_m$ is a
basis of local sections of $\hat{\Delta}$ there exist local functions
$c_0^1,\dots,c_0^m \in C^\infty(M_0)$ such that
\begin{equation}
\hat{X}(\gamma(t)) - \{ \cdot, H_0\}_0(\gamma(t)) = c_0^1(\gamma(t))\hat{X}_1(\gamma(t)) + \ldots + c_0^m(\gamma(t)) \hat{X}_m(\gamma(t)) \in \hat{\Delta}(\gamma(t)), \; \forall t.
\end{equation}
Consider the derivation 
\begin{equation}
\hat{Y} = \{ \cdot, H_0\}_0  + c_0^1\hat{X}_1 + \ldots + c_0^m\hat{X}_m
\end{equation}
on $C^\infty(M_0)$. This is the projection of a local vector field $Y
= \{ \cdot,H \} + c^1 X_1 + \ldots + c^m X_m$ on $M$, where $c^j \in
C^\infty(M)^G$ are such that $c_0^j \circ \pi = c^j\vert_{P^{-1}(0)}$
and $\{X_1,\ldots,X_m\}$ denotes the projectable local basis of
$\Delta$. Since the flow of this vector field commutes with the
$G$-action, it preserves the submanifold $P^{-1}(0) \cap M_{(K)}$
(following the argument above Remark \ref{rem:functions}). It follows
that the flow corresponding to the integral curve $\gamma(t)$
preserves the pieces $(M_0)_{(K)}$ and therefore $\gamma(t)$ restricts
to a smooth curve $\bar{\gamma}(t)$ on $(M_0)_{(K)}$. The vector field
$Y$ restricts to a vector field $\hat{Y}^\prime$ on $(M_0)_{(K)}$. By
construction it follows that
\begin{equation}
\hat{Y}^\prime = \{ \cdot, H_{(K)} \}_{(K)} + \bar{c}_0^1 \hat{X}_1^\prime + \ldots +\bar{c}_0^m \hat{X}_m^\prime,
\end{equation}
where $\bar{c}_0^j = c_0^j\vert_{(M_0)_{(K)}}$ and $\hat{X}_j^\prime$
is the restriction of $\hat{X}_j$ to $(M_0)_{(K)}$, $j=1,\ldots,m$.
The curve $\bar{\gamma}(t)$ is an integral curve of $\hat{Y}^\prime$.
Indeed,
\begin{equation}
\frac{d}{dt} \bar{f}_0 (\bar{\gamma}(t)) = \frac{d}{dt} f_0(\gamma(t)) = \hat{X}[f_0](\gamma(t)) = \hat{Y}[f_0](\gamma(t)) = \hat{Y}^\prime[\bar{f}_0](\bar{\gamma}(t)), \; \forall \bar{f}_0 \in W^\infty((M_0)_{(K)}),
\end{equation}
and since $W^\infty((M_0)_{(K)})$ is dense in $C^\infty((M_0)_{(K)})$
the result follows. Furthermore, it is clear that
\begin{equation}
\hat{Y}^\prime (\bar{\gamma}(t)) - J_{(K)}(\bar{\gamma}(t))dH_{(K)}(\bar{\gamma}(t)) = \hat{Y}^\prime (\bar{\gamma}(t)) - \{ \cdot,H_{(K)} \}_{(K)}(\bar{\gamma}(t)) \in \Gamma_{loc}\left( (\hat{\Delta})_{(K)}\right)(\bar{\gamma}(t)), \; \forall t,
\end{equation}
and
\begin{equation}
\hat{Z}^\prime [H_{(K)}](\bar{\gamma}(t)) = 0, \; \forall t, \; \forall \hat{Z}^\prime \in \Gamma_{loc}\left((\hat{\Delta})_{(K)}\right),
\end{equation}
which means that $\bar{\gamma}(t)$ is a solution of the regular
reduced Hamiltonian system $((M_0)_{(K)},\D_{(K)},H_{(K)})$. In
conclusion we have proved the following:

\begin{proposition}
  Consider the conditions in Proposition \ref{prop:pieces}. The
  regular reduced implicit Hamiltonian systems
  $((M_0)_{(K)},\D_{(K)},H_{(K)})$ are exactly the restrictions of the
  singular reduced implicit Hamiltonian system $(M_0,\D_0,H_0)$ to the
  pieces $(M_0)_{(K)}$. A solution $\gamma(t)$, with $\gamma(0) \in
  (M_0)_{(K)}$, of $(M_0,\D_0,H_0)$ preserves the piece $(M_0)_{(K)}$
  and restricts to a solution $\bar{\gamma}(t)$ of
  $((M_0)_{(K)},\D_{(K)},H_{(K)})$.
\end{proposition}

Finally, remark that since the pieces $(M_0)_{(K)}$ are smooth
manifolds, each implicit Hamiltonian system
$((M_0)_{(K)},\D_{(K)},H_{(K)})$ can be written as a set of
differential and algebraic equations (DAE). The singular reduced
implicit Hamiltonian system $(M_0,\D_0,H_0)$ can thus be written as a
collection of DAEs, one (set of differential and algebraic equations)
on each piece.

\section{Conclusions}\label{sec:concl}
In this paper we studied the singular reduction of implicit
Hamiltonian systems admitting a symmetry Lie group with a
corresponding equivariant momentum map. The results extend the
singular reduction theory known for explicit symplectic or Poisson
Hamiltonian systems \cite{ACG91,BL97,CB97,CS91,CSn91,OR98,OR03,SL91}. The
main result is a purely topological description of the reduced
implicit Hamiltonian system using the definition of a topological
Dirac structure. In particular, the reduced space is not assumed to be
a smooth manifold. The dynamics corresponding to this system are
defined and it is shown that the projectable solutions of the
unreduced system project to solutions of the reduced system. If the
symmetry Lie group acts freely and properly and the value of the
momentum map is regular, then the singular reduced implicit
Hamiltonian system equals the regular reduced implicit Hamiltonian
system as described in \cite{B00,BvdS01}. Finally, under certain
conditions, the singular reduced space can be decomposed into a set of
smooth manifolds called pieces. It is shown that the singular reduced
implicit Hamiltonian system restricts to regular reduced implicit
Hamiltonian systems on all the pieces.

\bigskip
 
\noindent\textbf{Acknowledgments.} This research was performed during
the time when the first author held a postdoctoral position at EPFL,
Switzerland. G.B. gratefully acknowledges the hospitality and the
financial support of this institution. Furthermore, G.B. also
acknowledges the financial support received from the European
sponsored project GeoPlex (IST-2001-34166, \texttt{www.geoplex.cc}).
T.S.R. was partially supported by the European Commission and the
Swiss Federal Government through funding for the Research Training
Network \emph{Mechanics and Symmetry in Europe} (MASIE) as well as the
Swiss National Science Foundation.


\end{document}